\documentclass[a4paper,fleqn]{article}

\usepackage[utf8]{inputenc}
\usepackage[english]{babel}

\usepackage{graphicx}

\renewcommand{\d}{\mathrm{d}}
\newcommand{\e}{\mathrm{e}}

\usepackage{csquotes}
\usepackage{xspace}

\usepackage{amsmath}
\DeclareMathOperator\arctanh{arctanh}

\usepackage{amsfonts}
\usepackage[usenames,dvipsnames]{xcolor}

\usepackage{subcaption}
\usepackage{subfiles}
\usepackage{wrapfig}
\usepackage{todonotes}

\usepackage[top=1in, bottom=1.25in, left=1.25in, right=1.25in]{geometry}

\begin{document}

\title{Derivation of Delay Equation Climate Models Using the Mori-Zwanzig Formalism}

\author{
Swinda K.J. Falkena$^{1,2}$, Courtney Quinn$^{3,4}$, Jan Sieber$^3$, Jason Frank$^{5,6}$
\\ and Henk A. Dijkstra$^{1,6}$}

\maketitle

\noindent{
$^{1}$Institute for Marine and Atmospheric Research Utrecht, Department of Physics, Utrecht University, Utrecht, The Netherlands\\
$^{2}$Mathematics of Planet Earth Program, University of Reading, Reading, UK\\
$^{3}$Department of Mathematics, University of Exeter, Exeter, UK\\
$^{4}$CSIRO Oceans and Atmosphere, Hobart, TAS, AU\\
$^{5}$Mathematical Institute, Utrecht University, Utrecht, The Netherlands\\
$^{6}$Centre for Complex Systems Studies, Faculty of Science, Utrecht University, Utrecht, The Netherlands\\
}

{\centering
\emph{Corresponding author}: {s.k.j.falkena@student.reading.ac.uk}\\[3mm]

\emph{Subject}: {Applied mathematics, Oceanography, Differential equations}\\[3mm]

}

\begin{abstract}
Models incorporating delay have been frequently used to understand climate variability phenomena, but often the delay is introduced through an ad-hoc physical reasoning, such as the propagation time of waves.  In this paper, the Mori-Zwanzig formalism is introduced as a way to systematically derive delay models from systems of partial differential equations and hence provides a better justification for using these delay-type  models. The Mori-Zwanzig technique gives a formal rewriting of the system using a  projection onto a set of resolved variables, where the rewritten system contains a memory term. The computation of this memory term requires solving the orthogonal dynamics equation, which represents the unresolved dynamics. For nonlinear systems, it is often not possible to obtain an analytical solution to the orthogonal dynamics and an approximate solution needs to be found. Here, we demonstrate the Mori-Zwanzig technique for a two-strip model of the El Ni\~no Southern Oscillation (ENSO) and explore methods to solve the orthogonal dynamics. The resulting nonlinear delay model contains an additional term compared to previously proposed ad-hoc conceptual models. This new term leads to a larger ENSO period, which is closer to that seen in observations. 
\end{abstract}


%
%





\section{Introduction}
\label{sec:intro}

To study climate variability and climate change, a hierarchy of models is currently used \cite{DijkstraBook}. At the low end of this hierarchy are conceptual climate models, which contain only the necessary features for specific phenomena to occur and thus represent the dominant physical processes. These  models are often used to study physical mechanisms in their purest forms, such as the causal chain behind a specific oscillation and how the period of this oscillation depends on the different processes involved. At the high end of the hierarchy, there are modern multi-process, multi-scale global climate models (GCMs) that aim to represent the total climate system in substantial detail. These models are used, for example, to make projections of future climate change as documented by assessments of the Intergovernmental Panel on Climate Change (IPCC). 

One special class of conceptual models consists of differential delay models, which are particularly useful for compactly representing physical phenomena. Compared to ordinary differential equation models, delay models can potentially convey more information -- they are infinite-dimensional dynamical systems -- but can still be formulated in terms of functions of a single variable. This allows for an easier mathematical treatment than would a partial differential equation model, while such a delay model can still represent complex physical processes. A useful way in which these models can be analysed is through bifurcation analysis. Such analysis allows for the distinction of different  dynamical regimes  and the dependence on a few parameters can be investigated relatively easy compared to GCMs. 

Keane et al. \cite{Keane2017} provided an overview of delay models used to describe climate processes. The two main areas in the climate system upon which delay models so far have focused are Energy Balance Models and models for the El Ni\~no Southern Oscillation (ENSO). There exist many positive and negative feedbacks in the climate system. Some of these feedbacks are delayed by a non-negligible amount of time, for example by transport through an ocean basin. The presence of such a delayed feedback can sometimes be determined from data \cite{Runge2014}. In models incorporating delay it is not necessary to resolve all of the processes involved in the feedback. A parametrization by the resulting delay time is sufficient. Usually delay models for climate  phenomena  are derived from more complex models by making strong assumptions about the  system, see e.g. \cite{Suarez1988, Battisti1988a}.  In many studies, the delay model is introduced in an ad-hoc way, usually through physical reasoning or semi-empirical indications. For instance, in many processes in climate, propagating waves play an important role and a wave-basin crossing time is used as a delay. 

A prominent case of delayed feedback through wave propagation is the ENSO variability in the Tropical Pacific.  During an El Ni\~no event, the sea surface temperature in the eastern part of the basin is warmer than usual, resulting in severe weather disturbances in countries on both sides of the Pacific Ocean. Its counterpart is La Ni\~na, when the sea surface temperature is colder than usual. These two events alternate with intermediate phases in between, resulting in an irregular oscillation with a period of four to seven years. One of the most successful models of ENSO is that by Zebiak and Cane \cite{Zebiak1987}. The view of the behaviour of ENSO in terms of normal modes resulted in the so-called delayed-oscillator mechanism of ENSO \cite{Jin1997a,Jin1997b}. The delay mechanism here is related to the propagation of equatorial Kelvin and off-equatorial Rossby waves. These waves take time to travel through the basin, resulting in a delayed arrival of a temperature anomaly. Already in 1988,  Suarez and Schopf proposed a delay model  for the sea surface temperature $T$ in the eastern Pacific Ocean of the form \cite{Suarez1988}:
\begin{equation}
\label{eq:ssm}
\frac{\d T}{\d t} = T(t) - T(t)^3 - \alpha T(t-\delta).
\end{equation}
Here $\delta$ is the delay time and $\alpha$ a parameter indicating the strength of the delayed feedback. Note that this model is scaled to contain as few parameters as possible. Other delay models based on the same mechanism have been proposed and studied, see e.g. \cite{Tziperman1998, Krauskopf2014}.

The delay model of Equation (\ref{eq:ssm}) gives oscillations with a period of two to three years for realistic values of $\alpha$ and $\delta$. This is on the short side with respect to the observed period of ENSO, indicating that some aspect is missing in the model. The physics behind the delay is well understood, justifying the linear terms in Equation (\ref{eq:ssm}). However, the nonlinear term in the model by Suarez and Schopf is proposed ad-hoc and no physical justification is given for it in their article. Battisti and Hirst provided some arguments for the form of the nonlinearity \cite{Battisti1988a}, but no thorough mathematical derivation has yet been provided.

The formal model reduction approach proposed by Mori \cite{Mori1965} and Zwanzig \cite{Zwanzig1973} potentially provides a way to place the derivation of delay models upon a stronger mathematical foundation. Using the Mori-Zwanzig approach, one formally reduces the dimension of a system of ordinary differential equations by projecting the dynamics onto a select subset of resolved dependent variables \cite{Chorin2002, Givon2004}. Closure is attained by replacing dependence on the unresolved variables by a memory integral and a term referred to as the noise term. This noise term only represents noise when the unresolved variables are chaotic on their own. Its statistics then are determined by the initial values of the unresolved variables. 

A common application of the formalism is stochastic modelling \cite{Gottwald2017}, for example for systems with large scale differences, which in certain cases can be reduced to Markovian systems \cite{Wouters2016}. For Markov chains there are several approaches to model reduction \cite{Beck2009}. Also for non-Markovian systems some results exist \cite{Darve2009}. Hamiltonian systems are another class of systems to which the formalism has been applied \cite{Chorin2002, Zhu2018a}. The component of the rewritten system that is focused on in this paper is the memory term. When the memory term is expressed as the convolution of a kernel function with the history of the system, it is in the form of a distributed delay. Under some approximations, this memory term can be simplified to a term with a discrete delay.

In this paper the Mori-Zwanzig formalism is applied to a spatially extended model of ENSO, which serves as a test case for the derivation of delay models using this formalism. In Section \ref{sec:mz} the theory behind the Mori-Zwanzig formalism is  shortly recapitulated. In Sections \ref{ssec:twostr} and \ref{ssec:nonlts} the formalism is applied to linear and nonlinear two-strip models of ENSO and new delay models are derived. The new nonlinear models are analyzed in Section \ref{ssec:delan} and the results are summarized and discussed in Section \ref{sec:disc}.

\section{Mori-Zwanzig Formalism}
\label{sec:mz}

The Mori-Zwanzig formalism reformulates a set of ordinary differential equations (ODEs) into a reduced system for the resolved variables that still retains all of the dynamics of the original system. In this section, we outline the formalism. We adopt the notation and follow the approach of Chorin \emph{et al.} \cite{Chorin2002} in the formulation of the reduction, which is based on the work by Mori (1965) \cite{Mori1965} and Zwanzig (1973) \cite{Zwanzig1973}. The theory discussed in this section holds for ODEs, whereas the system we will consider is a partial differential equation (PDE). The formalism is not straightforward to apply to general systems of PDEs, which are infinite-dimensional. Problems arising when treating PDEs are discussed in Section \ref{ssec:ENSOmz}.

We start by presenting a simple linear example with constant coefficients, to illustrate the idea behind the formalism. Consider the following linear system of ODEs for $\phi = (\hat{\phi}, \tilde{\phi}): \mathbb{R} \rightarrow \mathbb{R}^n$ continuously differentiable:
\begin{equation}
\label{eq:linmz1}
\frac{\d }{\d t}
\begin{pmatrix}
\hat{\phi} \\
\tilde{\phi}
\end{pmatrix} = 
\begin{pmatrix}
A_{11} & A_{12} \\
A_{21} & A_{22}
\end{pmatrix}
\begin{pmatrix}
\hat{\phi} \\
\tilde{\phi}
\end{pmatrix}
\end{equation}
Here, we call $\hat{\phi}\in\mathbb{R}^m$ the \emph{resolved} variables and $\tilde{\phi}\in\mathbb{R}^{n-m}$ the \emph{unresolved} variables.  We have $A_{11} \in \mathbb{R}^{m\times m}$, $A_{12} \in \mathbb{R}^{m\times (n-m)}$, $A_{21} \in \mathbb{R}^{(n-m)\times m}$ and $A_{22} \in \mathbb{R}^{(n-m)\times(n-m)}$ and initial conditions $\phi(0) = (\hat{x}, \tilde{x})$. The goal is to derive an equation for the resolved variables $\hat{\phi}$ only. In this example, this can be done by solving the equation for $\tilde{\phi}$ using variation of constants and substituting the result into the equation for $\hat{\phi}$:
\begin{equation}
\label{eq:linmz2}
\frac{\d}{\d t} \hat{\phi}(t) = A_{11} \hat{\phi}(t) + A_{12} \e^{A_{22}t}\tilde{x} + \int_0^t A_{12} \e^{A_{22}(t-s)} A_{21} \hat{\phi}(s) \d s.
\end{equation}
The system in Equation (\ref{eq:linmz1}) has been reduced to one equation for the resolved variables $\hat{\phi}$, with the only dependence on $\tilde{\phi}$ being the  initial condition $\tilde{x}$. This reduced system is equivalent to the full system and exhibits the same behaviour. The first term on the right-hand side is referred to as the Markovian term, the second as the noise term (due to possible uncertainty in the initial condition) and the last as the memory term.

For the generalization of this idea to nonlinear systems, consider the following $n$-dimensional system of ODEs:
\begin{equation}
\label{eq:odemzth}
\frac{\d}{\d t}\phi(t) = R(\phi(t)), \quad \phi(0) = x,
\end{equation}
where $\phi(t)\in\mathbb{R}^n$ is a continuously differentiable function of $t\in\mathbb{R}_+$, $x\in\mathbb{R}^n$ denotes the initial condition, and $R:\mathbb{R}^n\rightarrow\mathbb{R}^n$ has components $R_i$. To every initial condition $x$ there corresponds a trajectory $\phi(t)=\phi(x,t)$, $\phi:\mathbb{R}^n \times \mathbb{R}_+ \rightarrow \mathbb{R}^n$, whose existence is assumed for all $t>0$. We consider the evolution of an observable $u(x,t) := g(\phi(x,t))$ along a solution of Equation (\ref{eq:odemzth}), where $g$ is defined on $\mathbb{R}^n$. The quantity $u(x,t)$, $u:\mathbb{R}^n \times \mathbb{R}_+ \rightarrow \mathbb{R}^n$, satisfies the PDE
\begin{equation}
\label{eq:liouville}
\frac{\partial}{\partial t} u(x,t) = \mathcal{L} u(x,t), \quad u(x,0) = g(x),
\end{equation}
where $[\mathcal{L}u](x) = \sum_{i=1}^n R_i(x)\partial_{x_{i}}u(x)$ is the generator associated with vector field $R$ of Equation (\ref{eq:odemzth}). This generator is called the Liouville operator \cite{morriss2013}. 

The goal, as for the linear system, is to construct a system of equations for a select subset of $m$ resolved variables $\hat{\phi}\in\mathbb{R}^m$. As before, the unresolved variables are denoted by $\tilde{\phi}\in\mathbb{R}^{n-m}$, such that $\phi = (\hat{\phi}, \tilde{\phi})$. To reduce the system from $n$ components to the desired $m$ components a projection operator $P:C(\mathbb{R}^n,\mathbb{R}^k)\rightarrow C(\mathbb{R}^m,\mathbb{R}^k)$ is needed. Here, $k$ is the dimension of an arbitrary function $f$ to which the projection is applied. Examples of projection operators include the conditional expectation (infinite-rank) \cite{Dominy2017} and the linear projection, defined by $[Pf](\hat{x}) = f(\hat{x},0) =: \hat{f}(\hat{x})$, which sets all unresolved variables to zero and retains only the resolved components. The complement of $P$ is denoted $Q=I-P$, where $I$ is the identity operator. 

We denote the solution of the linear PDE (\ref{eq:liouville}) as $u(x,t) = [\e^{t\mathcal{L}} g](x)$, where $\e^{t\mathcal{L}}$ is also referred to as the evolution operator \cite{Darve2009}.  In particular for $g(x) = x_i$ we find $\phi_i(x,t) = \e^{t\mathcal{L}}x_i$.  Using this notation, Equation (\ref{eq:liouville}) can be written as
\begin{equation}
\label{eq:lioq}
\frac{\partial}{\partial t} [\e^{t\mathcal{L}}g](x) = [\e^{t\mathcal{L}} \mathcal{L}g] (x) = [\e^{t\mathcal{L}} P \mathcal{L}g] (x) + [\e^{t\mathcal{L}} Q \mathcal{L} g](x).
\end{equation}
Note that $\mathcal{L}$ and $\e^{t\mathcal{L}}$ commute. We now consider the second term in the right-hand side of this equation, $[\e^{t\mathcal{L}} Q \mathcal{L} g](x)$. This component gives the evolution of the unresolved variables. The Dyson formula \cite{morriss2013},
\begin{equation}
\e^{t(A+B)} = \e^{tA} + \int_0^t \e^{(t-s)(A+B)}B \e^{sA} \d s,
\end{equation}
applied for $A = Q\mathcal{L}$ and $B = P\mathcal{L}$ gives
\begin{equation}
[\e^{t\mathcal{L}} Q\mathcal{L}g](x) = [\e^{tQ\mathcal{L}}Q\mathcal{L}g](x)   + 
\int_0^t [\e^{(t-s)\mathcal{L}}P\mathcal{L} \e^{sQ\mathcal{L}} Q\mathcal{L}g](x) \d s.
\end{equation}
Substitution into Equation (\ref{eq:lioq}) yields
\begin{equation}
\label{eq:prelan}
\frac{\partial}{\partial t}  [\e^{t\mathcal{L}}g](x) = [\e^{t\mathcal{L}}P\mathcal{L}g](x) + [\e^{tQ\mathcal{L}}Q\mathcal{L}g](x)   + 
\int_0^t [\e^{(t-s)\mathcal{L}}P\mathcal{L} \e^{sQ\mathcal{L}} Q\mathcal{L}g](x) \d s.
\end{equation}
In particular for $g(x) = x_i$ we find $[\e^{t\mathcal{L}}g](x)=\phi_i(x,t)$ and Equation \ref{eq:prelan} becomes the generalized Langevin equation:
\begin{equation}
\label{eq:lan}
\frac{\partial}{\partial t} \phi_i(x,t) = R_i(\hat{\phi}(x,t))  + F_i(x,t) + \int_0^t K_i(\hat{\phi}(x,t-s),s) \d s,
\end{equation}
where we use the shorthand notation $R_i(\hat{\phi}(x,t))= R_i([\hat{\phi}(x,t),0]) =[PR_i](\phi(x,t))$ and
\begin{equation}
\label{eq:defnoimem}
F_i(x,t) = [\e^{tQ\mathcal{L}} Q \mathcal{L} g](x), \qquad K_i(\hat{x},t) = [P\mathcal{L} F_i](x,t),
\end{equation}
with notation $K_i(\hat{x},t) = K_i([\hat{x},0],t)$. Note that $F_i(x,t)$ is the solution to the orthogonal dynamics equation:
\begin{equation}
\label{eq:ortd}
\frac{\partial}{\partial t}F_i(x,t) = Q \mathcal{L} F_i(x,t), \qquad F_i(x,0) = Q \mathcal{L} x_i.
\end{equation}
In general it is not known whether this initial value system is well-posed, but in specific cases (approximate) solutions can be obtained. The three terms on the right-hand side of the Langevin equation (\ref{eq:lan}) are called the Markovian term $R_i(\hat{\phi}(x,t))$, the noise term $F_i(x,t)$ and the memory term, defined as the integral over the memory integrand $K_i(\hat{\phi}(x,t-s),s)$. This memory integrand consists of a memory kernel applied to the resolved variables. Note that the integral over $K_i(\hat{\phi}(x,s'),t-s')$ is equal to the memory term in Equation \eqref{eq:lan} after a change of variables $s'=t-s$. From this point onward we will use the integral in $s'$, dropping the prime for conciseness. For linear systems solving the orthogonal dynamics system simplifies to the case of Equation \eqref{eq:linmz1} and the Langevin equation is equivalent to the result obtained by using variation of constants.

The generalized Langevin system (\ref{eq:lan}) is still exact, but it is not necessarily simpler. If solving the orthogonal dynamics equation (\ref{eq:ortd}) is as difficult as solving the full system, there is no use in applying the formalism. The applicability thus depends on the particular system and whether a suitable projection exists. Such a projection would yield an orthogonal dynamics system that is relatively straightforward to solve or approximate in a good way. Applications to slow-fast \cite{Wouters2016}, Markovian \cite{Beck2009}, non-Markovian \cite{Darve2009}, and Hamiltonian systems \cite{Chorin2002}, as well as systems with an orthogonal basis of eigenfunctions \cite{Szalai2014}, have been considered in the literature \cite{Givon2004}. For non-Hamiltonian systems of PDEs, which are considered here, less is known. The standard approach when considering PDE systems is to expand in different modes of the system \cite{Szalai2014}. More recently an expansion based on the Faber series has been proposed \cite{Zhu2018b}.

The main challenge when applying the Mori-Zwanzig formalism is the solution of the orthogonal dynamics system (\ref{eq:ortd}). The choice of projection operator $P$ is an important factor in determining the form and complexity of the orthogonal dynamics equation. The projection should be chosen such that the orthogonal dynamics system is stable, meaning you need to retain stabilizing factors in the unresolved dynamics, and less complex than the original system. Alternatively, one may approximate the orthogonal dynamics equation by a less complex system. An example is the pseudo-orthogonal dynamics approximation derived by Gouasmi \emph{et al.} \cite{Gouasmi2017}, which is discussed in Appendix A. This approximation is applied in Sections \ref{ssec:ENSOmz} and \ref{ssec:ENSOnlmz}. Preferably, the orthogonal dynamics system decays at a faster rate than the full system. For linear constant coefficient systems this means that the largest eigenvalue of the orthogonal dynamics system is smaller than that of the full system. In that case one may justify neglecting the noise term $F_i(x,t)$ in the Langevin equation.

\section{Linear ENSO  Model}
\label{ssec:twostr}

In this section we apply the Mori-Zwanzig formalism to a linear model of ENSO variability. The model we study is a system of PDEs in one space variable describing the dynamics on two strips, one at the equator and one at higher latitude, in the Pacific Ocean. The Mori-Zwanzig formalism and the use of characteristics reduce this model to a linear delay equation similar to the model by Suarez and Schopf \cite{Suarez1988}.

\subsection{Model Formulation}
\label{ssec:formlin}

ENSO is a coupled ocean-atmosphere phenomenon, where the variations in the sea surface temperature (SST) induce wind stress anomalies which drive ocean circulation changes affecting the SST. The interaction between the wind stress and the ocean is described in the two-strip ocean model derived and studied by Jin \cite{Jin1997a,Jin1997b}. This two-strip model is derived from the dimensionless shallow water equations in a normalized equatorial basin \cite{DijkstraBook}. Assuming a parabolic dependence of the thermocline on latitude near the equator simplifies the shallow water equations to a system of equations for the thermocline depth at the equator ($h_e$) and at some latitude $y_n$ between 5$^\circ$N and 15$^\circ$N ($h_n$). Since the model used is scaled, the basin is of length one with $x=0$ being the western boundary and $x=1$ the eastern boundary of the Pacific ocean at the equator. The details of this scaling can be found in Appendix B.

The thermocline depth at both latitudes responds to a wind forcing whose strength depends on the ocean temperature at the equator. To describe this coupling we use a simplified version of the Gill atmosphere model \cite{DijkstraBook}. An equation for the SST perturbations $T_e$ at the equator, see e.g. \cite{Dijkstra1995}, completes the system. The resulting two-strip model describing the dynamics of ENSO is \cite{DijkstraBook}
\begin{subequations}
	\label{eq:twostr} 
	\begin{align}
	(\partial_t + \epsilon_0) (h_e - h_n) + \partial_x h_e & = \mu g(x) T_e(x_E,t),  \label{eq:twostr1} \\
	(\partial_t + \epsilon_0) h_n - \frac{1}{y_n^2}\partial_x h_n & = -\mu \frac{\theta}{y_n^2} g(x) T_e(x_E,t), \label{eq:twostr2} \\
	\partial_t T_e + c_T T_e - c_h h_e & = 0, \label{eq:twostr3} 
	\end{align}
\end{subequations}
for $x\in[0,1]$. The boundary conditions are 
\begin{equation}
h_e(0,t) = r_W h_n(0,t), \qquad h_n(1,t) = r_E h_e(1,t).
\end{equation}
Here $\epsilon_0$ is a linear damping coefficient, $\mu$ a coupling coefficient for the wind forcing and $\theta$ an order one coefficient representing the difference in the effect of wind-stress between the equator and higher latitudes. In the equation for SST the coefficient $c_T$ represents local damping and $c_h$ represents the effect of thermocline depth on temperature through background upwelling. Both $c_T$ and $c_h$ can depend on time, space or any of the variables ($h_c$, $h_n$, $T_e$), which can result in nonlinearities. Furthermore, $r_W$ and $r_E$ are a measure of the mass flux at the western and eastern boundaries respectively. The function $g(x)$ gives the pattern of the wind forcing in the zonal direction. The wind forcing depends on the SST anomaly at $x=x_E$, in the east of the basin. Because the SST only feeds back into the thermocline equations for $x=x_E$, it is sufficient for solving the system to only consider $T_e$ at that location. The $\partial_x$-terms represent the advection of anomalies in the thermocline by Kelvin waves ($h_e$) and Rossby waves ($h_n$). 

The left-hand side of Equations \eqref{eq:twostr1} and \eqref{eq:twostr2} can be decoupled by introducing a new variable $h_c = h_e - \frac{1}{1+y_n^2}h_n$ and keeping $h_n$. The new variable $h_c$ is dominated by the thermocline depth at the equator, but includes some influence of the higher latitudes as well. In the new variables the equations are
\begin{subequations}
	\label{eq:tworw}
	\begin{align}
	\partial_t h_c + \epsilon_0 h_c + \partial_x h_c &= \mu \Big(1-\frac{\theta}{1+y_n^2}\Big) g(x) T_e(x_E,t), \label{eq:tworw1}\\
	\partial_t h_n + \epsilon_0 h_n - \frac{1}{y_n^2}\partial_x h_n &= -\mu \frac{\theta}{y_n^2} g(x) T_e(x_E,t), \label{eq:tworw2}\\
	\partial_t T_e + c_T T_e - c_h \Big(h_c + \frac{1}{1+y_n^2} h_n\Big) &= 0, \label{eq:tworw3}
	\end{align}
\end{subequations}
with boundary conditions
\begin{equation}
\label{eq:bctwostr}
h_c(0,t) = \Big(r_W - \frac{1}{1+y_n^2}\Big) h_n(0,t), \qquad r_E h_c(1,t) = \Big(1-\frac{r_E}{1+y_n^2}\Big) h_n(1,t).
\end{equation}
If $r_E=0$, then $h_n(1,t)=0$, meaning no reflection occurs at the eastern boundary.

The solution to the homogeneous equations for thermocline depth in the rewritten system of Equations (\ref{eq:tworw}), that is without wind forcing ($\mu=0$), can be expanded in eigenmodes for rates $\sigma_k$: 
\begin{equation}
\label{eq:eigtwo}
h_c^0(x,t) = H_c \e^{\sigma_k t} \e^{-(\sigma_k+\epsilon_0)x}, \qquad
h_n^0(x,t) = H_n \e^{\sigma_k t} \e^{(\sigma_k+\epsilon_0)y_n^2 x},
\end{equation}
with
\begin{equation}
\sigma_k = -\epsilon_0 + \frac{1}{1+y_n^2} \Big(\ln \Big( \frac{r_Er_W(1+y_n^2) - r_E}{(1+y_n^2) - r_E} \Big) + 2\pi i k\Big), \qquad k\in\mathbb{N}.
\end{equation}
The boundary conditions \eqref{eq:bctwostr} imply that $H_c = \Big(r_W - \frac{1}{1+y_n^2}\Big) H_n$, where $H_n$ is arbitrary. For $r_E=0$ the solution of the homogeneous system is trivial after one round trip: $h_c^0 = h_n^0 = 0$. The solutions in Equation (\ref{eq:eigtwo}) are the eigensolutions of the two-strip model. Note that these eigensolutions are not orthogonal, meaning they are not convenient to use as a basis on which can be projected.

\subsection{Mori-Zwanzig Formalism}
\label{ssec:ENSOmz}

Starting from the rewritten version of the two-strip model in Equation (\ref{eq:tworw}), the goal is to derive a delay equation describing ENSO using the Mori-Zwanzig formalism. The resulting model is expected to be similar to that by Suarez and Schopf \cite{Suarez1988} as given in Equation (\ref{eq:ssm}). The Mori-Zwanzig formalism works on ODEs, while the two-strip model is a system of PDEs. A route that is often taken in such situations is to expand the general solution in a basis of eigensolutions of the PDE \cite{Darve2009,Givon2004,Szalai2014}. Truncating this expansion and projecting along the adjoint eigensolutions reduces the PDE to a system of ODEs. However, for the two-strip model this is impractical, since the eigensolutions may be degenerate. Other projection methods, such as the use of Fourier exponentials or orthogonal polynomials as a basis \cite{Shen2011}, will converge very slowly (thus, requiring high truncation order) due to incompatible boundary conditions, making them unsuitable for analytical computation.

Here we treat $\partial_x$ as an operator and consider the system as ODEs on the space of continuous functions (in $x$). The result of the Mori-Zwanzig formalism will contain terms of the form $\e^{t\partial_x}$, which are solution operators of scalar wave equations. In this section we consider the linear version of the two-strip model, meaning the coefficients $c_T$ and $c_h$ are allowed to depend on space ($x$) and possibly time ($t$), but not on any of the variables $h_c$, $h_n$ or $T_e$. Specifically, we permit dependencies $c_T(x)$ and $c_h(x,\epsilon t)$, where $\epsilon$ represents the option for a slowly varying background effect of thermocline anomalies on SST. Note that by considering only linear terms, we cannot expect to find the cubic term in the model by Suarez and Schopf in Equation \eqref{eq:ssm}. Extension to a nonlinear model is discussed in Section \ref{ssec:nonlts}.

The first step in applying the Mori-Zwanzig formalism is to identify the Liouville operator following its definition in Section \ref{sec:mz}:
\begin{equation}
\label{eq:liouts}
\begin{split}
\mathcal{L} &= \Big(-\big(\epsilon_0 + \partial_x\big) h_c(x,0) + \mu \Big(1-\frac{\theta}{1+y_n^2}\Big) g(x) T_e(x_E,0)\Big)\partial_{h_c}  \\
& \quad + \Big(-\big(\epsilon_0 -\frac{1}{y_n^2} \partial_x\big) h_n(x,0) - \mu \frac{\theta}{y_n^2} g(x) T_e(x_E,0)\Big)\partial_{h_n} \\
& \quad + \Big(-c_T(x) T_e(x,0) + c_h(x,\epsilon t)\big(h_c(x,0) + \frac{1}{1+y_n^2} h_n(x,0)\big) \Big)\partial_{T_e}.
\end{split}
\end{equation}
The second step is to choose a projection operator $P$. The model by Suarez and Schopf is a delay equation for the temperature at the equator in the east of the basin. Therefore, we choose the equatorial temperature $T_e$ as the resolved variable and use the linear projection as defined in Section \ref{sec:mz},
\begin{equation}
\label{eq:plin}
P(f(T_e,h_c,h_n)) = f(T_e,0,0) =: \hat{f}(T_e),
\end{equation}
with $P:C(\mathbb{R}^3,\mathbb{R})\rightarrow C(\mathbb{R},\mathbb{R})$, reducing the number of dependent variables from three to one. This projection is infinte-rank as $T_e$ is still a function of $x$. Applying the formalism thus results in an equation for only the resolved variable $T_e$.

Now we have the information needed to apply the formalism and find the Langevin equation \eqref{eq:lan} for the linear two-strip model. The different terms are computed and discussed separately. Firstly the Markovian term is computed:
\begin{equation}
\label{eq:mark}
[\e^{t\mathcal{L}}P\mathcal{L}T_e](x,0)
= -c_T(x) T_e(x,t).
\end{equation}
As expected, this is the right-hand side dependence on the resolved variable $T_e$ in Equation \eqref{eq:tworw3}, which gives the evolution of the resolved variable in the model.

To compute the noise and memory term we need to solve the orthogonal dynamics equation. One problem that can arise with infinite-rank projections is that it has not been generally proven that $\e^{tQ\mathcal{L}}$ is well-posed. The projection we consider has infinite rank, however since $P\mathcal{L}Q$ is bounded, the operator $\e^{tQ\mathcal{L}}$ can be bounded by $M\e^{\omega t}$ for some $\omega$ \cite{Zhu2018a} (note that $g(x)$ in Equation (\ref{eq:liouts}) is empirically derived and will be bounded - see Section \ref{ssec:dellin}). This implies the contribution of the memory term is bounded and $\e^{tQ\mathcal{L}}$ is well-posed.

Because the operator $Q \mathcal{L}$ is linear in this case, we can use the method of Gouasmi \emph{et al.} \cite{Gouasmi2017} to simplify the orthogonal dynamics system in Equation \eqref{eq:ortd}. For an explanation of this simplification the reader is referred to Appendix A. The resulting orthogonal dynamics are advection equations with complete boundary conditions. The equations for the thermocline depth are
\begin{equation}
\label{eq:ortdyn_enso}
\begin{split}
\partial_t h_c^Q(x,t) &= -(\epsilon_0+\partial_x)h_c^Q(x,t), \\
\partial_t h_n^Q(x,t) &= -(\epsilon_0-\frac{1}{y_n^2}\partial_x)h_n^Q(x,t),
\end{split}
\end{equation}
with boundary conditions \eqref{eq:bctwostr} for $h_c^Q$ and $h_n^Q$.  The equation for $T_e^Q$ decouples and is not needed for computation of the noise or memory term. Here $Q$ is used to denote the variables in the orthogonal dynamics equation. These two equations are independent of each other and have exponential solutions. The solutions are the same as those of the homogeneous two-strip model as given in Equation \eqref{eq:eigtwo}. Having solved the orthogonal dynamics equation, the noise term can be computed following Equation \eqref{eq:defnoimem}:
\begin{equation}
\label{eq:ensolinnoise}
\begin{split}
F_{T_e}(x,t) &= c_h(x,\epsilon t) \big(h_c^Q(x,t) + \frac{1}{1+y_n^2} h_n^Q(x,t)\big) \\
&= c_h(x,\epsilon t) \big(e^{-(\epsilon_0 + \partial_x)t} h_c(x,0) + \frac{1}{1+y_n^2} e^{-(\epsilon_0 - \frac{1}{y_n^2}\partial_x)t} h_n(x,0)\big).
\end{split}
\end{equation}
Here the solutions of the orthogonal dynamics system \eqref{eq:ortdyn_enso} have been substituted to find the final expression.

The last component of the Langevin equation that needs to be computed is the memory term. To do this, first the memory integrand is computed following Equation (\ref{eq:defnoimem}):
\begin{equation}
\label{eq:ensolinmemory}
\begin{split}
K_{T_e}(T_e(x,0),t) &= c_h(x,\epsilon t) \Big(\mu \Big(1-\frac{\theta}{1+y_n^2}\Big) \e^{-(\epsilon_0 + \partial_x)t}  - \mu \frac{\theta}{y_n^2} \frac{1}{1+y_n^2} \e^{-(\epsilon_0 - \frac{1}{y_n^2}\partial_x)t} \Big) \\
& \qquad \cdot g(x) T_e(x_E,0).
\end{split}
\end{equation}
Here we exploit that $P$ and $\mathcal{L}$ commute with $\e^{-(\epsilon_0 + \partial_x)t}$ and $\e^{-(\epsilon_0 - \frac{1}{y_n^2}\partial_x)t}$. This can be verified by comparing with the result of applying variation of constants (Equation \eqref{eq:linmz2}). As expected $K_{T_e}$ is linear in $T_e$. The memory term in the Langevin equation is found by substituting the computed memory integrand into the integral.

The resulting Langevin equation is obtained by substituting the computed Markovian \eqref{eq:mark}, noise \eqref{eq:ensolinnoise} and memory term \eqref{eq:ensolinmemory} into Equation \eqref{eq:lan}. We find the following equation for the temperature at the equator:
\begin{subequations}
	\label{eq:mzresenso}
	\begin{align}
	\frac{\d T_e}{\d t}(x,t) &= -c_T(x) T_e(x,t) \label{enso:lin:markov}\\
	& \quad + c_h(x, \epsilon t)\Big( \e^{-(\epsilon_0 + \partial_x)t} h_c(x,0) + \frac{1}{1+y_n^2} \e^{-(\epsilon_0 - \frac{1}{y_n^2}\partial_x)t} h_n(x,0)\Big) \label{enso:lin:noise}\\
	& \quad + \int_0^t c_h(x,\epsilon (t-s)) \Big(\mu \big(1-\frac{\theta}{1+y_n^2}\big) \e^{-(\epsilon_0 + \partial_x)(t-s)} \nonumber\\
	& \quad \qquad - \mu \frac{\theta}{y_n^2}\frac{1}{1+y_n^2} \e^{-(\epsilon_0 - \frac{1}{y_n^2}\partial_x)(t-s)} \Big) g(x) T_e(x_E,s) \d s. \label{enso:lin:memory}
	\end{align}
\end{subequations}
Note the change of variables discussed in Section \ref{sec:mz} has been applied to arrive at the memory term \eqref{enso:lin:memory}. The partial derivative to $x$ is still present in the exponential terms. In the noise term \eqref{enso:lin:noise} this is not an issue, since the terms are exactly the solutions to the homogeneous system as given in Equation (\ref{eq:eigtwo}). In the memory term \eqref{enso:lin:memory} it is less clear how to evaluate these exponential operators. In the next section these terms are simplified using the method of characteristics.

\subsection{Evaluation along Characteristics}
\label{ssec:char}

In this section we focus on the exponential $\partial_x$-terms in the memory integral \eqref{enso:lin:memory} using the method of characteristics \cite{PDEBook}. Both components in the memory kernel are of the form $e^{-(\epsilon_0 + c\partial_x)(t-s)}  f(x,s)$, for either $c=1$ or $c=-1/y_n^2$. This expression is the solution to the PDE
\begin{equation}
\partial_t f + c \partial_x f = -\epsilon_0 f,
\end{equation}
with initial conditions given at $t=s$. The characteristic curves of this equation are $x-x_0 = c(t-t_0)$, along which $f(x,t)$ is constant apart from damping caused by the $\epsilon_0$-term. This gives
\begin{equation}
\label{eq:chargen}
\e^{-(\epsilon_0 + c\partial_x)(t-s)}  f(x,s) = \e^{-\epsilon_0(t-s)} f(c(t-s)+x_s,s),
\end{equation}
where $x_s$ is the location at time $s$. This expression is valid as long as the argument of $f$ lies in the domain. Otherwise boundary effects need to be considered.

\begin{figure}[h]
	\centering
	\includegraphics[width=0.4\textwidth]{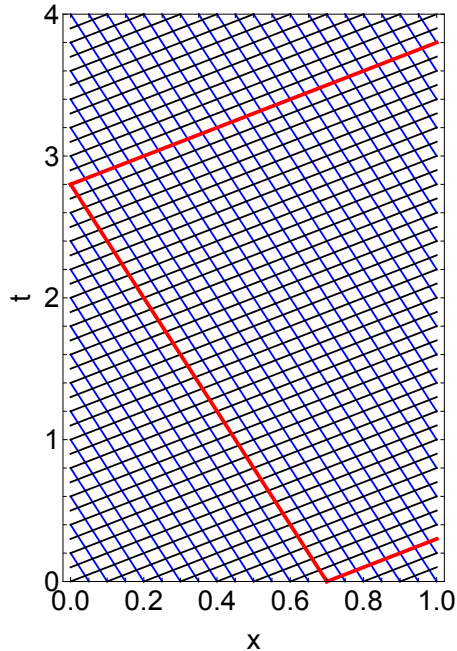}
	\caption{\label{fig:char}The characteristics of $\partial_t f + \partial_x f = -\epsilon_0 f$ (black) and $\partial_t f - \frac{1}{y_n^2}\partial_x f = -\epsilon_0 f$ (blue). In red the path of a signal following the characteristics is shown until it reaches the eastern boundary.}
\end{figure}

In the memory integral the two exponential terms have different corresponding characteristics. They represent the eastward traveling equatorial Kelvin waves for $c=1$ and the westward traveling Rossby waves for $c=-1/y_n^2$. Note that the Rossby waves take longer to cross the basin than the Kelvin waves as $y_n>1$. In Figure \ref{fig:char} the characteristics of the memory term are shown. The red line shows what happens to a signal emitted from $x=0.7$ at time zero until it arrives at the eastern boundary.

Since the domain of the two-strip model is bounded ($x\in[0,1]$), the effect of these boundaries on the temperature signal needs to be discussed. Following the boundary conditions in Equation \eqref{eq:bctwostr}, we can express the equation for $T_e$ at $x_b=0,1$ in terms of either $h_n(x_b,t)$ for the westward Rossby wave, or $h_c(x_b,t)$ for the eastward Kelvin wave. If we consider $r_E$ and $r_W$ in Equation \eqref{eq:bctwostr} to be nonzero (i.e.~allowing for energy transfer between strips at the boundaries) then reflection of the characteristics must be considered.  The fraction between the incoming and outgoing wave, corrected for their respective effects on $T_e$ (Equation \eqref{eq:tworw3}), determine the reflection coefficients for the western and eastern boundary which are given by, respectively, 
\begin{equation}
\label{eq:twobc}
A_{rW} = r_W(1+y_n^2) - 1, \qquad \text{and} \qquad A_{rE} = \Big(\frac{1+y_n^2}{r_E}-1\Big)^{-1}.
\end{equation}

The characteristics of the system are used to get an expression for the memory term. The goal is to find a result for the temperature in the east of equatorial basin ($x=x_E$), where the model by Suarez and Schopf is defined \cite{Suarez1988}. This also is the location of the temperature on which the wind forcing depends. Looking at the signal at one location allows for the following of characteristics from a source given by $g(x)$ to that one location. The memory integral at $x=x_E$ is of the form
\begin{equation}
\label{eq:charxE}
\begin{split}
&\int_0^t c_h(x_E,\epsilon (t-s)) \Big(\mu \big(1-\frac{\theta}{1+y_n^2}\big) \cdot \Big[\e^{-(\epsilon_0 + \partial_x)(t-s)}  g(x)\Big]_{x_E} T_e^E(s) \\
& \qquad - \mu \frac{\theta}{y_n^2}\frac{1}{1+y_n^2} \cdot \Big[\e^{-(\epsilon_0 - \frac{1}{y_n^2}\partial_x)(t-s)} g(x)\Big]_{x_E} T_e^E(s) \Big) \d s.
\end{split}
\end{equation}
The first term in the memory integral \eqref{enso:lin:memory} gives the waves that are traveling westward at $t=s$. Setting $x_E=1$ for the characteristics, these waves need a time $t=1-x$ to arrive at the eastern boundary starting from an initial location $x$. Assuming reflection takes place at the eastern and western boundary, the signal arrives a second time after $t=1-x+(y_n^2+1)$. The signal keeps reflecting through the basin, arriving at the eastern boundary after times $t_k=1-x+k(y_n^2+1)$ for $k=0,1,2,...$. At each reflection the wave loses energy by a factor $A_{rE}$ at the eastern boundary and a factor $A_{rW}$ at the western boundary.

The result of these reflections through the basin is found by combining Equation \eqref{eq:chargen} and the above discussion in Equation \eqref{enso:lin:memory}. The resulting expression for the first part of the memory term at the eastern boundary is
\begin{equation}
\begin{split}
& \sum_{k=0}^{K_{max}(t)} \int_{t-(1+k(y_n^2+1))}^{t-k(y_n^2+1)} c_h(x_E,\epsilon (t-s)) \mu \big(1-\frac{\theta}{1+y_n^2}\big) \e^{-\epsilon_0(t-s)} \\
&\qquad \cdot (A_{rE}A_{rW})^k g(1+k(y_n^2+1)-(t-s)) T_e^E(s) \d s,
\end{split}
\end{equation}
where $K_{max}(t) = \lfloor \frac{t-1}{y_n^2+1} \rfloor$ for $t\geq 1$ is the number of reflections that have occurred by time $t$. Note that there are time intervals for which this term has no effect at the eastern boundary, since it only represents one of the two characteristics. To get this part of the memory integral in a form which shows more of the delay behaviour, a change of coordinates can be applied. Let $x=1+k(y_n^2+1)-(t-s)$, for which $\frac{\d x}{\d s}=1$. Changing coordinates from $s$ to $x$, yields a memory integral for $T_e^E(t-(1+k(y_n^2+1)-x))$. This shows that the memory term contains a component that depends on past states of the resolved variable $T_e^E$. How strong the effect is at a certain time depends on the function $g(x)$, which gives the spatial distribution of the wind forcing.

The result for the waves that first travel towards the western boundary is achieved in a similar way. These waves need times $t_k=y_n^2 x +1+k(1+y_n^2)$ for $k=0,1,2,...$ to arrive at the eastern boundary. Now $x = -\frac{1}{y_n^2}(1+k(y_n^2+1)-(t-s))$ is used for the change of coordinates. Going through the same steps as before, the total memory integral becomes
\begin{equation}
\label{eq:memch}
\begin{split}
& \sum_{k=0}^{K_{max}(t)} (A_{rE}A_{rW})^k \mu \int_{0}^{1} g(x) \e^{-\epsilon_0k(y_n^2+1)}\\
&\cdot \Big( \big(1-\frac{\theta}{1+y_n^2}\big) c_h(x_E,\epsilon (1+k(y_n^2+1)-x)) \e^{-\epsilon_0 (1-x)} T_e^E(t-(1+k(y_n^2+1)-x)) \\
& - \frac{\theta}{y_n^2}\frac{A_{rW}}{1+y_n^2} c_h(x_E,\epsilon (1+k(y_n^2+1)+y_n^2x)) \e^{-\epsilon_0 (1+ y_n^2 x)} T_e^E(t-(1+k(y_n^2+1)+y_n^2x)) \Big) \d x.
\end{split}
\end{equation}
This expression shows that there are multiple delays present in the two-strip model. The exact form of the delay (distributed or discrete) and delay times are determined by the spatial pattern of the wind forcing $g(x)$.

\subsection{Delay Model}
\label{ssec:dellin}

Using the results from the two previous sections, a delay equation for the evolution of the SST in the east of the basin can be obtained. We consider Equation (\ref{eq:mzresenso}) with the expression for the memory term in Equation (\ref{eq:memch}) at $x=x_E$. The effect of the components of the memory term decreases for higher $k$ by energy loss at reflection. To simplify this expression we make two assumptions.

First, we assume that there is no reflection at the eastern boundary, meaning $r_E=0$ and thus $A_{rE}=0$. The only two components in the sum of the memory term (Equation \eqref{eq:memch}) that remain with this assumption are the components for $k=0$. As noted in Section \ref{ssec:formlin}, the homogeneous solution in that case is identically zero after finite time, such that the noise term \eqref{enso:lin:noise} vanishes.

The function $g(x)$ in the memory integral is unspecified until now. It determines the form of the memory kernel by setting the pattern of the wind forcing. More specifically, $g(x)$ indicates where the effect of the wind is strong and weak. Following Jin \cite{Jin1997a} the wind dominantly has an effect near the centre of the basin. Away from this location the effect is small, meaning that the wind forcing acts quite locally. We approximate this local effect of the wind forcing by $g(x) = A_0 \delta_{x_w}(x)$, a delta function of height $A_0$ at $x=x_w$. This delta function leaves only the effect of that one location on the integral, such that Equation (\ref{eq:mzresenso}) considered at $x=x_E$ simplifies to
\begin{equation}
\label{eq:ensolin1}
\begin{split}
\frac{\d T_e^E}{\d t} &= -c_T(x_E) T_e^E(t) + \mu A_0 \Big( \big(1-\frac{\theta}{1+y_n^2}\big) c_h(x_E,\epsilon (1-x_w)) \e^{-\epsilon_0 (1-x_w)} T_e^E(t-(1-x_w)) \\
& \qquad - \frac{\theta}{y_n^2}\frac{A_{rW}}{1+y_n^2} c_h(x_E,\epsilon (1+y_n^2x_w)) \e^{-\epsilon_0 (1+ y_n^2 x_w)} T_e^E(t-(1+y_n^2x_w)) \Big).
\end{split}
\end{equation}
This is a linear equation with discrete delay for the temperature. Because we assumed $c_h$ and $c_T$ to be independent of $h_{c}$, $h_n$ and $T_e$ (see \ref{ssec:ENSOmz}) no nonlinearity is found.

Equation \eqref{eq:ensolin1} does not yet resemble the model by Suarez and Schopf from Equation (\ref{eq:ssm}) in the linear terms. Instead of one, there are two delay times present. However, $1-x_w \ll 1+y_n^2x_w$ for realistic parameters (see Appendix B) indicating that this effect can be considered to be immediate. Thus it is assumed that $T_e^E(t-(1-x_w))\approx T_e^E(t)$.  This approximation yields the final linear delay model:
\begin{equation}
\label{eq:dellin}
\begin{split}
\frac{\d T_e^E}{\d t} &= c_S T_e^E(t) - c_L T_e^E(t-d),
\end{split}
\end{equation}
where
\begin{equation}
\begin{split}
c_S &= \mu A_0 \big(1-\frac{\theta}{1+y_n^2}\big) c_h(x_E,\epsilon (1-x_w)) \e^{-\epsilon_0 (1-x_w)} - c_T(x_E), \\
c_L &= \mu A_0 \frac{\theta}{y_n^2}\frac{A_{rW}}{1+y_n^2} c_h(x_E,\epsilon (1+y_n^2x_w)) \e^{-\epsilon_0 (1+ y_n^2 x_w)}, \\
d &= 1+y_n^2x_w.
\end{split}
\end{equation}
This model (after rescaling) gives the linear part of the model by Suarez and Schopf in Equation \eqref{eq:ssm} \cite{Suarez1988}. The delay is due to the propagation  of Rossby waves caused by a wind forcing which depends on the temperature near the eastern boundary. These waves travel to the western boundary, where they reflect in the form of Kelvin waves. This delay model does not yet account for the nonlinearity in the model by Suarez and Schopf, without which no stable oscillation will occur in the model. To get a more realistic result a nonlinear version of the two-strip model is considered as the starting point for applying the Mori-Zwanzig formalism in the next section.

\section{Nonlinear ENSO  Model}
\label{ssec:nonlts}

In this section we start by deriving a nonlinear variation of the two-strip model describing the ENSO dynamics. This nonlinearity allows for a more realistic modelling of ENSO. The resulting nonlinear two-strip model is studied using the Mori-Zwanzig formalism and variation of constants.

\subsection{Model Formulation}
\label{ssec:formnonlin}

In the two-strip model \eqref{eq:tworw}, as considered in Section \ref{ssec:twostr}, the thermocline feedback coefficient $c_h$ in the temperature equation did not depend on any of the variables $T_e$, $h_e$ or $h_n$. We now introduce a realistic state dependence into $c_h$, resulting in a nonlinear PDE system. The expression for $c_h$ is given by \cite{Dijkstra1995}:
\begin{equation}
\label{eq:twoch}
c_h = f_h(x) \frac{\d T_s}{\d h},
\end{equation}
where $f_h(x)$ is the background wind forcing and $T_s$ the subsurface temperature at the equator as a function of thermocline depth $h$. For the parametrization of the subsurface temperature the result by Hao \emph{et al.} is used \cite{Hao1993}:
\begin{equation}
\label{eq:Tsch}
T_s(h) = T_{s0} + (T_0 - T_{s0}) \tanh\Big( \frac{h+h_0}{H^*} \Big).
\end{equation}
Here $T_0$ is the ocean equilibrium temperature in absence of dynamics, $h_0$ an offset value for the thermocline and $T_{s0}$ the temperature at $h=-h_0$. The parameter $H^*$ determines the steepness of the transition when $h$ passes through $-h_0$. The derivative of $T_s$ to $h$, which is needed to get an expression for $c_h$ \eqref{eq:twoch}, is proportional to $1-\Big(\frac{T_s-T_{s0}}{T_0-T_{s0}}\Big)^2$.

\begin{figure}[h]
	\centering
	\includegraphics[width=.6\textwidth]{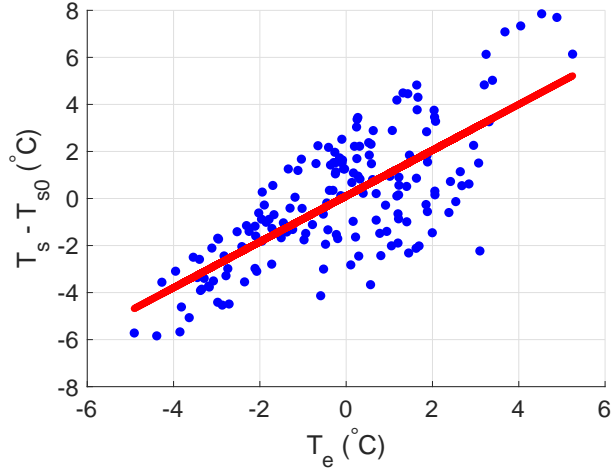}
	\caption{\label{fig:cortemp}Temperature data from a buoy at the eastern side of the Pacific Ocean (110$^\circ$W) for measurements at the surface and the subsurface (depth of 40m). Shown are anomalies of the SST versus anomalies of the subsurface temperature (average 20.5 $^\circ$C). The red line is the best linear fit through the data and has a slope of $0.97\pm0.07$. Data is taken from the Global Tropical Moored Buoy Array Project Office of NOAA/PMEL.}
\end{figure}

We assume that $T_e$ is proportional to $T_s-T_{s0}$, meaning perturbations in the equatorial sea surface temperature are proportional to perturbations in the subsurface temperature. Using this assumption, Equation (\ref{eq:twoch}) has the form
\begin{equation}
\label{eq:chnlfull}
c_h(x, T_e) = f_h(x)\frac{T_0-T_{s0}}{H^*}\Big( 1 - \Big( \frac{c_{se} T_e}{T_0-T_{s0}}\Big)^2 \Big),
\end{equation}
where $c_{se}$ is the proportionality constant. This will introduce a cubic nonlinearity in the temperature equation of the two-strip model, which was not present in the linear version considered in Section \ref{ssec:twostr}. To check the validity of the assumption, buoy data of the equatorial Pacific Ocean is considered. Since the model only contains feedback of the temperature in the east of the basin, it is sufficient to consider only buoys in the eastern part of the Pacific Ocean. There are ten locations in the equatorial Pacific where buoy data is available. To avoid coastal boundary layer effects the second most eastern buoy, which is located at 110 degrees west, is chosen. In Figure \ref{fig:cortemp} anomalies of the SST versus those of the subsurface temperature (40 m) are shown for this buoy. The correlation between the two datasets is 0.83, indicating there is a strong relation. The slope between the two temperature anomalies is $c_{se}\approx 1$. We note that also at other depths correlation is strong. However, the slope between the temperature anomalies decreases slightly with depth.

Using Equation (\ref{eq:chnlfull}) the nonlinear two-strip model (in the form of Equation \eqref{eq:tworw}) becomes
\begin{equation}
\label{eq:twonl}
\begin{split}
\partial_t h_c + \epsilon_0 h_c + \partial_x h_c &= \mu \Big(1-\frac{\theta}{1+y_n^2}\Big) g(x) T_e(x_E,t), \\
\partial_t h_n + \epsilon_0 h_n - \frac{1}{y_n^2}\partial_x h_n &= -\mu \frac{\theta}{y_n^2} g(x) T_e(x_E,t), \\
\partial_t T_e + c_T(x) T_e - c_h^*(x)(1-\beta T_e^2) \Big(h_c + \frac{1}{1+y_n^2} h_n\Big) &= 0,
\end{split}
\end{equation}
where $c_h^*(x)=f_h(x)\frac{T_0-T_{s0}}{H^*}$ and $\beta=\Big( \frac{1}{T_0-T_{s0}}\Big)^2$. We emphasize that this result is in principle only valid in the eastern part of the basin, which is sufficient for the model presented here. At other locations the correlation between anomalies in SST and subsurface temperatures is quite strong as well, but the proportionality constant is different. In the following we discuss two ways to apply the Mori-Zwanzig formalism to the nonlinear system of Equation \eqref{eq:twonl}. The first method is based on variation of constants, the second on an approximation to the Mori-Zwanzig formalism.

\subsection{Variation of Constants}
\label{ssec:voc}

Since the equations for $h_c$ and $h_n$ are still linear, we may apply the Mori-Zwanzig formalism in its elementary form \eqref{eq:linmz2} based on the variation of constants. This yields the following equation for $T_e$:
\begin{equation}
\begin{split}
\frac{\d T_e}{\d t} &=  - c_T(x) T_e(x,t) + c_h^*(x)(1-\beta T_e^2(x,t)) \\
&\qquad \cdot \Big(\e^{-(\epsilon_0+\partial_x)t} h_c(x,0) + \int_0^t \e^{-(\epsilon_0+\partial_x)(t-s)} \mu \Big( 1- \frac{\theta}{1+y_n^2}\Big) g(x) T_e(x_E,s) \d s  \\
&\qquad + \frac{1}{1+y_n^2} \Big( \e^{-(\epsilon_0-\frac{1}{y_n^2}\partial_x)t} h_n(x,0) - \int_0^t \e^{-(\epsilon_0-\frac{1}{y_n^2}\partial_x)(t-s)} \mu \frac{\theta}{y_n^2} g(x) T_e(x_E,s) \d s \Big)\Big).
\end{split}
\end{equation}
The terms within the integrals are the same as those in the memory integral of Equation \eqref{enso:lin:memory}. Therefore, we can use the evaluation along characteristics in Section \ref{ssec:char} to rewrite them in the form of Equation \eqref{eq:memch}.

Then, under the same assumptions as made in Section \ref{ssec:char}, that is no reflection at the eastern boundary and a localized wind forcing, we may apply equivalent simplifications to obtain a delay equation for $T_e^E$. Assuming, as before, that the short delay is instantaneous, yields a nonlinear delay model for the temperature in the east of the basin, which is an exact reduction of the nonlinear two-strip PDE \eqref{eq:twonl}: 
\begin{equation}
\label{eq:delnlalt}
\begin{split}
\frac{\d T_e^E}{\d t} &= (c_S^*-c_T(x_E)) T_e^E(t) - c_L^* T_e^E(t-d) - \beta c_S^* T_e^{E}(t)^3 + \beta c_L^* T_e^{E}(t)^2 T_e^E(t-d),
\end{split}
\end{equation}
where
\begin{equation}
\begin{split}
c_S^* &= \mu A_0 \big(1-\frac{\theta}{1+y_n^2}\big) c_h^*(x_E) \e^{-\epsilon_0 (1-x_w)}, \\
c_L^* &= \mu A_0 \frac{\theta}{y_n^2}\frac{A_{rW}}{1+y_n^2} c_h^*(x_E) \e^{-\epsilon_0 (1+ y_n^2 x_w)}, \\
d &= 1+y_n^2x_w.
\end{split}
\end{equation}
This is a nonlinear equation with discrete delay for the temperature at the equator, including two cubic terms. The difference with the model by Suarez and Schopf in Equation (\ref{eq:ssm}) is a fourth term, which is proportional to $T_e^{E}(t)^2 T_e^E(t-d)$. The effect of this additional term on the dynamics of the delay model is studied in Section \ref{ssec:delan}. Before discussing the behaviour of this extended delay model, first the application of the Mori-Zwanzig formalism in the form of the Langevin equation \eqref{eq:lan} to the nonlinear model in Equation \eqref{eq:twonl} is considered.

\subsection{Mori-Zwanzig Formalism}
\label{ssec:ENSOnlmz}

The Mori-Zwanzig formalism is valid for both linear and nonlinear equations. The challenge when considering nonlinear equations arises in solving the orthogonal dynamics equation. We use the same linear projection as for the linear model \eqref{eq:plin}. Here, it could be tempting to use a Taylor expansion to find a formal expression for the noise term. However, because the nonlinear operator $\mathcal{L}$ is unbounded, this does not apply \cite{Kato1995}. An alternative is the Faber expansion \cite{Zhu2018b}, which yields a numerically computable solution where the expansion is truncated after a certain number of terms. However, we aim at finding an analytic expression for the resulting scalar equation for $T_e^E$ as found in the linear case \eqref{eq:dellin}. For this purpose the truncation approximation is not suitable. Another option is to approximate the orthogonal dynamics equation by the pseudo-orthogonal dynamics (POD) equation as derived by Gouasmi \emph{et al.} \cite{Gouasmi2017}. The conditions for this approximation are not met for the nonlinear two-strip model, but it can be used as a first estimate.

Here we will apply this POD approximation to derive an approximate solution. This derivation is given in Appendix A. The resulting equation is simplified by using the method of characteristics, assuming no reflection takes place at the eastern boundary and considering a localized wind forcing, just as was done in Sections \ref{ssec:char} and \ref{ssec:dellin}. Considering the short delay as being instantaneous, the resulting nonlinear delay equation for the temperature in the east of the basin is
\begin{equation}
\label{eq:delnlmz}
\begin{split}
\frac{\d T_e^E}{\d t} &= (c_S^*-c_T(x_E)) T_e^E(t) - c_L^* T_e^E(t-d) - \beta c_S^* T_e^{E}(t)^3 + \beta c_L^* T_e^{E}(t-d)^3,
\end{split}
\end{equation}
where
\begin{equation}
\begin{split}
c_S^* &= \mu A_0 \big(1-\frac{\theta}{1+y_n^2}\big) c_h^*(x_E) \e^{-\epsilon_0 (1-x_w)}, \\
c_L^* &= \mu A_0 \frac{\theta}{y_n^2}\frac{A_{rW}}{1+y_n^2} c_h^*(x_E) \e^{-\epsilon_0 (1+ y_n^2 x_w)}, \\
d &= 1+y_n^2x_w.
\end{split}
\end{equation}
This equation is almost the same as Equation \eqref{eq:delnlalt}, which was derived by applying variation of constants. The difference is the fourth term. Since Equation (\ref{eq:delnlalt}) is exact, this is due to errors introduced by the POD approximation.

\section{Analysis of the Delay Models}
\label{ssec:delan}

Section \ref{ssec:twostr} shows that the presence of a delay follows already from the analysis of the two-strip model at the linear level, resulting in Equation \eqref{eq:ensolin1}. The nonlinearity due to temperature dependence of the thermocline feedback enters the resulting delay model in different terms depending on the reduction method (Section \ref{ssec:nonlts}). Equation \eqref{eq:delnlalt} is derived by applying variation of constants to the thermocline equations and is therefore exact. Equation \eqref{eq:delnlmz} is derived using the POD approximation of Gouasmi \emph{et al.} \cite{Gouasmi2017} (Appendix A) and contains some approximation error. Both models contain an extra term compared to the model proposed by Suarez and Schopf as given in Equation (\ref{eq:ssm}). In the rest of this section we will refer to Equation \eqref{eq:delnlalt} as the VoC model, Equation \eqref{eq:delnlmz} as the MZ model and Equation \eqref{eq:ssm} as the S\&S model.

Before studying the qualitative behaviour of the delay models \eqref{eq:delnlalt} and \eqref{eq:delnlmz} in more detail, we scale temperature and time to reduce the number of parameters. In this section the sub- and superscripts of temperature are omitted, such that $T$ is written for $T_e^E$. Time is scaled by $\tilde{t} = (c_S^*-c_T(x_E)) t$ and temperature by $\tilde{T} = \sqrt{\frac{\beta c_S^*}{c_S^*-c_T(x_E)}}T$, where $\tilde{t}$ and $\tilde{T}$ are the scaled quantities. The scaled equation for the VoC model of Equation (\ref{eq:delnlalt}) (omitting the tildes for simplicity) is
\begin{equation}
\label{eq:delscvoc}
\frac{\d T}{\d t} = T(t)-T(t)^3 - \alpha T(t-\delta) \big( 1 - \gamma T(t)^2 \big),
\end{equation}
where 
\begin{equation}
\alpha = \frac{c_L^*}{c_S^*-c_T(x_E)}, \qquad \gamma = \frac{c_S^*-c_T(x_E)}{c_S^*}, \qquad \delta = (c_S^*-c_T(x_E))d.
\end{equation}
For $\gamma=0$ Equation \eqref{eq:delscvoc} reduces to the S\&S model \eqref{eq:ssm}. Note that the scaled parameters do not depend on the strength of the nonlinearity in the thermocline feedback $\beta$. Only the temperature scale depends on $\beta$. Since $c_S^*>c_T(x_E)$ all scaled parameters are positive and $\gamma<1$. The scaling for the MZ model \eqref{eq:delnlmz}, obtained by the approximate Mori-Zwanzig formalism, is the same, resulting in:
\begin{equation}
\label{eq:delscmz}
\frac{\d T}{\d t} = T(t)-T(t)^3 - \alpha T(t-\delta) \big( 1 - \gamma T(t-\delta)^2 \big).
\end{equation}

\begin{figure}[h]
	\centering
	\includegraphics[width=.6\textwidth]{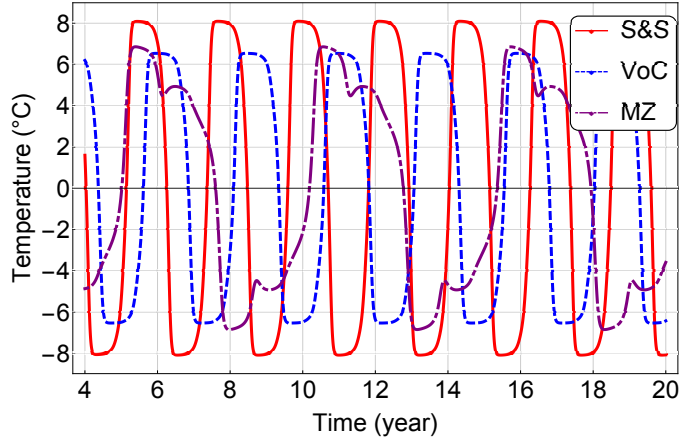}
	\caption{\label{fig:osc}Model simulations of delay models for ENSO: S\&S model (red), VoC model (blue) and MZ model (purple). The used parameters are: $\alpha=0.93$, $\gamma=0.49$ and $\delta=4.8$.}
\end{figure}

Figure \ref{fig:osc} shows time profiles for the three different nonlinear models for realistic parameters. The profile of the MZ model has a very different shape compared to the other two models. Both newly derived models have a longer period than the S\&S model, as well as a smaller amplitude. These two aspects are closer to measurements of ENSO. The main focus in this section is on the similarities and differences between the S\&S model by Suarez and Schopf and the VoC model derived by applying variation of constants.

\subsection{Bifurcation Analysis}
\label{ssec:bif}

Following the analysis by Suarez and Schopf \cite{Suarez1988}, we first determine the steady states of Equations \eqref{eq:delscvoc} and \eqref{eq:delscmz}. Depending on the value of $\alpha$, the models have either one or three equilibria:
\begin{equation}
T_{00} = 0, \qquad T_{0\pm} = \pm \sqrt{\frac{1 - \alpha}{1 - \alpha \gamma}} \quad \text{for } \alpha\notin\{1,1/\gamma\}.
\end{equation}
The trivial equilibrium, $T_{00} = 0$, undergoes a pitchfork bifurcation at $\alpha=1$, with two equilibria emerging for $\alpha<1$.  For $\alpha<1$, $T_{00}$ is always unstable, while for $\alpha>1$ the stability depends on the delay. 

\begin{figure}[h!]
	\centering
	\includegraphics[width=.8\textwidth]{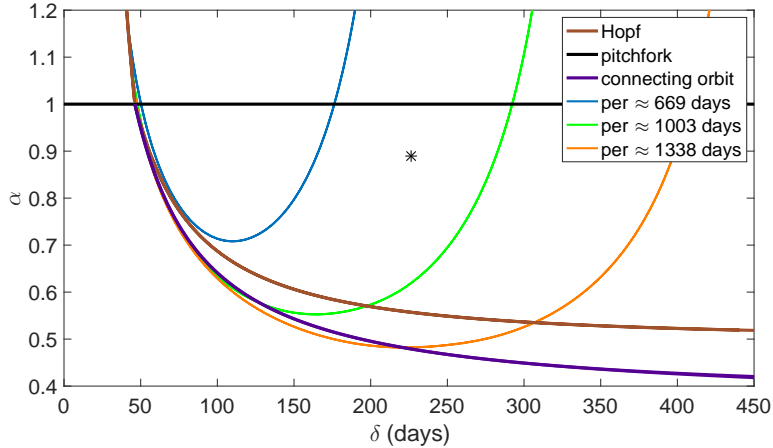}
	\caption{\label{fig:2param_SS}Bifurcation diagram in the $\alpha$-$\delta$-plane for the Suarez and Schopf \cite{Suarez1988} model. The black star shows the values used to compute Figure \ref{fig:osc}.}
\end{figure}

Using the MATLAB toolbox DDE-BIFTOOL \cite{DDEbiftool1,DDEbiftool2,DDEbiftool3} we compute two-parameter bifurcation diagrams in $\delta$ and $\alpha$, for each of the three nonlinear models.  Figure \ref{fig:2param_SS} shows the result for the original S\&S model (note that $\gamma=0$ makes the VoC and MZ models equivalent to the S\&S model).  For $\alpha>1$ and increasing $\delta$, $T_{00}$ undergoes a supercritical Hopf bifurcation with a family of stable periodic orbits emerging.  This large-amplitude symmetric periodic orbit corresponds to the ENSO behaviour. This family of ENSO periodic orbits extends to the region $\alpha<1$, where it emerges from a connecting orbit. The area between the connecting orbit and the Hopf curve for $\alpha<1$ is a multistable region, with the attractors being the ENSO periodic orbit and two stable equilibria. We also show level curves for the period of the ENSO orbit throughout the $\alpha$-$\delta$-plane. The values of the delay $\delta$ and the period length have been converted to the dimensional quantities for better comparison to observations (see Appendix B for scaling values).

\begin{figure}[h!]
	\centering
	\begin{subfigure}{.48\textwidth}
		\centering
		\includegraphics[width= \textwidth]{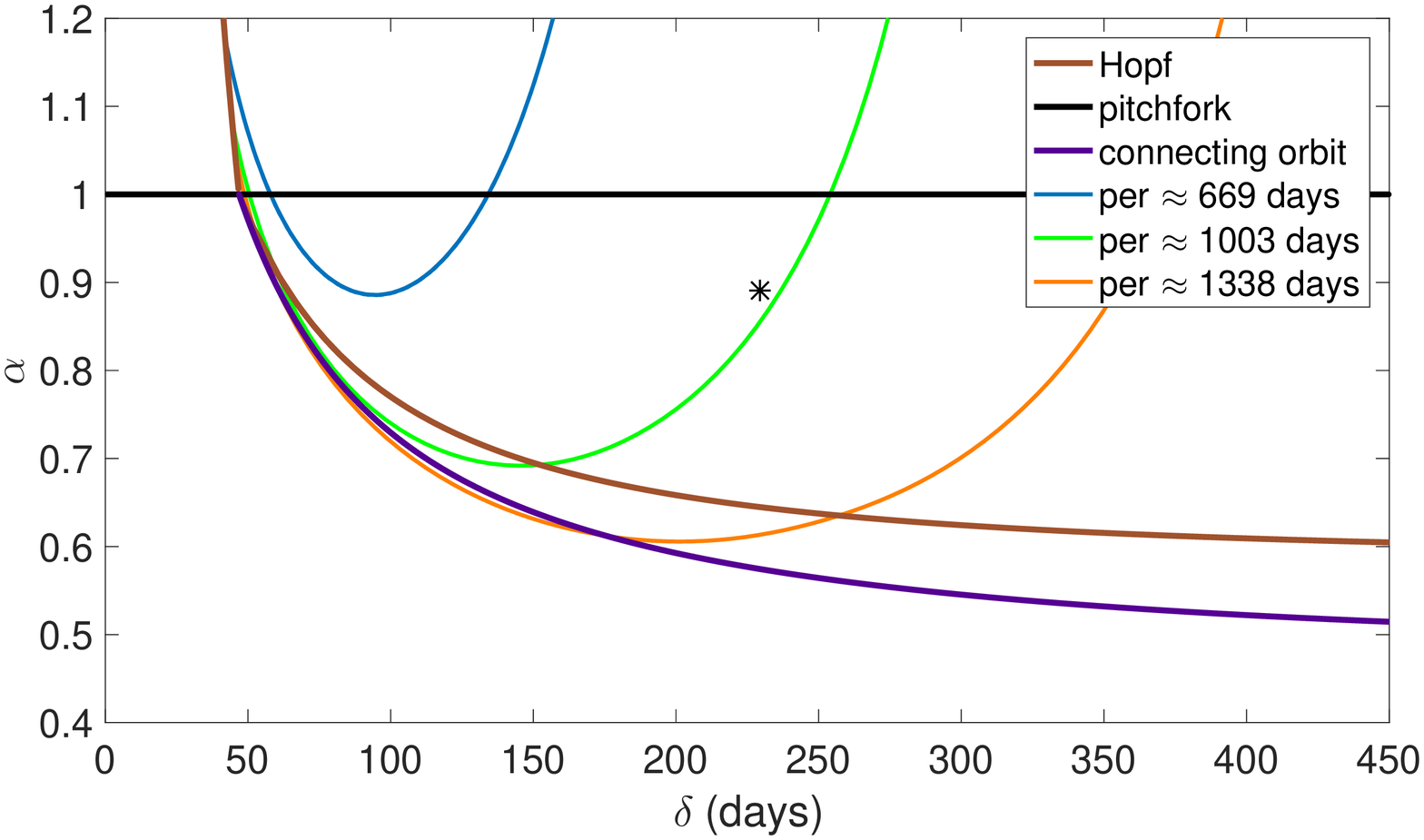}
		\caption{VoC \label{fig:bifgamVoC}}
	\end{subfigure}
	\begin{subfigure}{.48\textwidth}
		\centering
		\includegraphics[width= \textwidth]{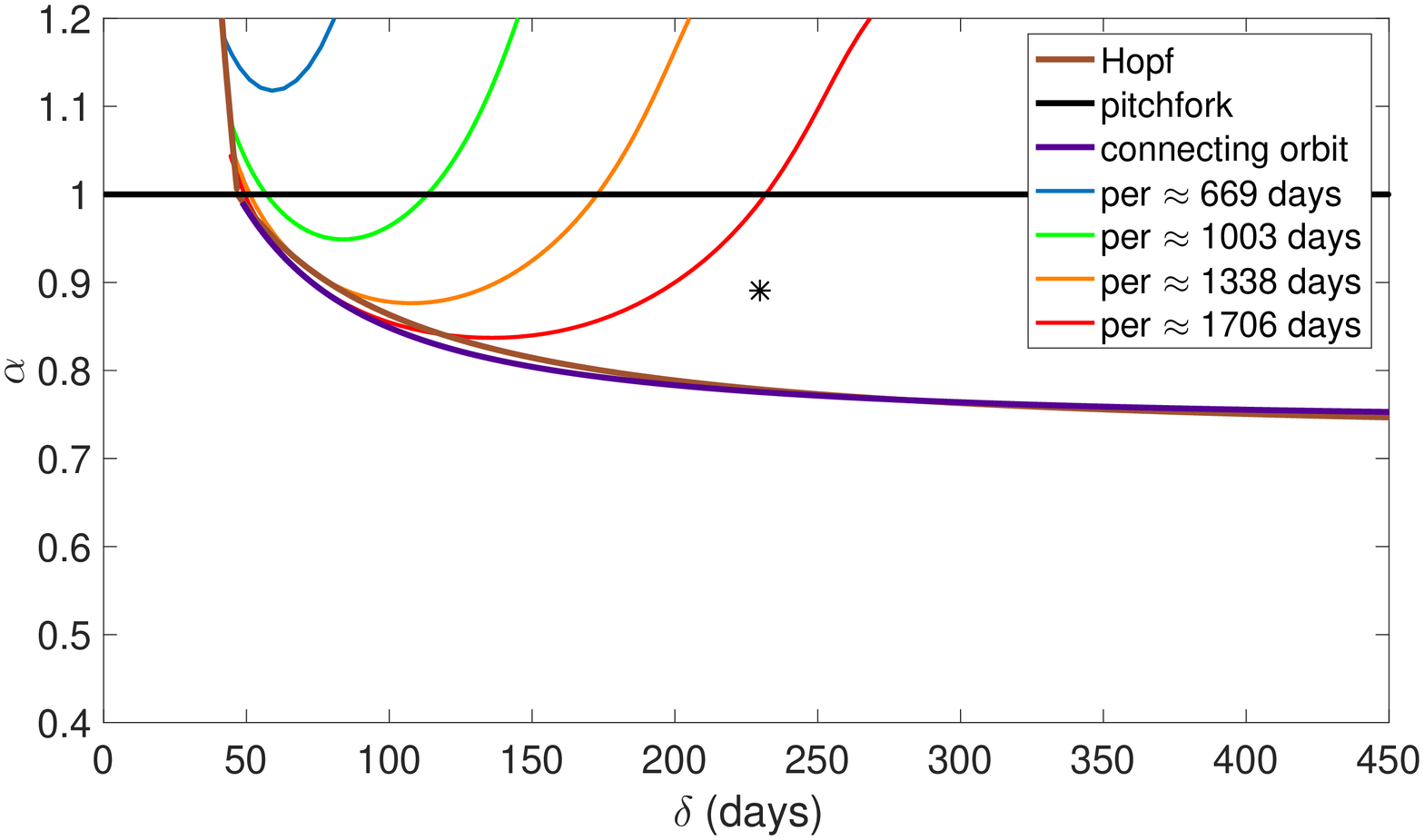}
		\caption{MZ}
	\end{subfigure}
	\caption{\label{fig:bifgamMZ}Bifurcation diagram in the $\alpha$-$\delta$-plane for the Variation of Constants (VoC) and Mori-Zwanzig (MZ) models. The black star shows the values used to compute Figure \ref{fig:osc}.}
\end{figure}

As shown in Section \ref{ssec:nonlts}, the extra terms involving $\gamma$ in the VoC and MZ model appear directly through the derivation from the ENSO PDE model.  Using the scaling in Appendix B, this would correspond to $\gamma=0.49$.  In Figure \ref{fig:bifgamMZ} we show the effect of incorporating non-zero $\gamma$ in the S\&S model through the same two-parameter bifurcation diagrams. We observe that in the VoC method (which is an exact reduction of the PDE), the region of parameter space for $\alpha<1$ where stable oscillations exist is slightly reduced. In the region between the Hopf curve and the connecting orbit curve in Figures \ref{fig:2param_SS} and \ref{fig:bifgamMZ} the ENSO oscillation coexists with two non-zero stable equilibria. However, the length of the period has increased (shown by the shift in the level curves),  as already suggested by Figure \ref{fig:osc}.

We can also get an idea of the dynamical consequences of the error in the POD approximation of Gouasmi \emph{et al.} \cite{Gouasmi2017} on the nonlinear ENSO model.  As mentioned in Section \ref{sec:mz}, the Mori-Zwanzig method is highly dependent on the ability to solve, or approximate to a good degree, the orthogonal dynamics equation.  The POD approximation of Gouasmi \emph{et al.} \cite{Gouasmi2017} (outlined in Appendix A) is applied in Section \ref{ssec:ENSOnlmz} to derive the MZ delay model.  We see through Figure \ref{fig:bifgamMZ} that this approximation introduces an error. The error causes an underestimation of the parameter region of stable oscillations for $\alpha<1$, and an overestimation of the period of oscillation.  In this particular case where one can find the exact projected equation through variation of constants, it is clear that the POD approximation is not sufficient and a different approximation should be proposed.  This is outside the scope of this study and should be considered in future work.

\subsection{Dependence of Period on Physical Parameters}
\label{ssec:pardep}

\begin{figure}
	\centering
	\includegraphics[width=.6\textwidth]{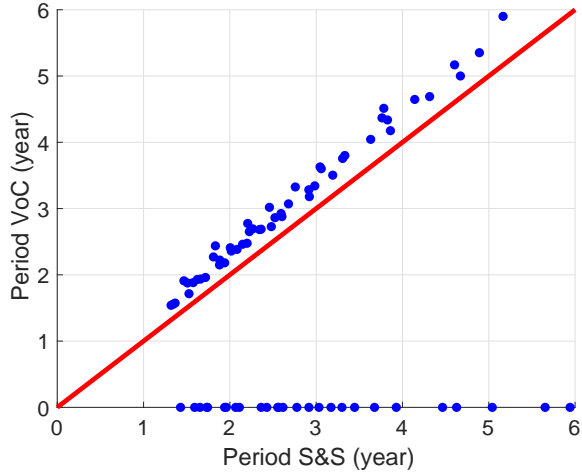}
	\caption{\label{fig:perch}The period of the VoC model versus the period according to the S\&S model. The line of equal period is shown in red.}
\end{figure}

\begin{figure}[h]
	\centering
	\begin{subfigure}{.32\textwidth}
		\centering
		\includegraphics[width= \textwidth]{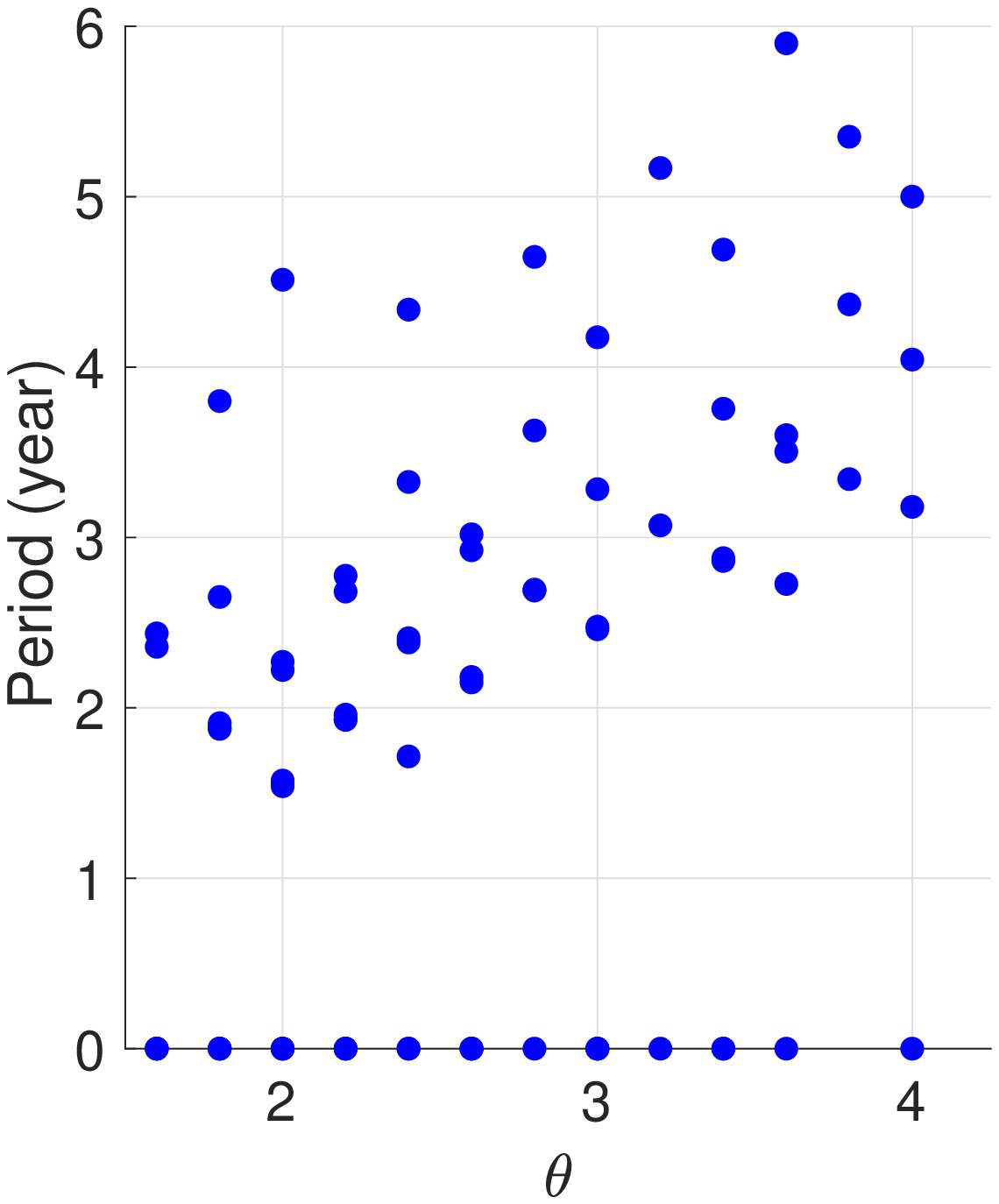}
		\caption{Changing $\theta$.}
	\end{subfigure}
	\begin{subfigure}{.32\textwidth}
		\centering
		\includegraphics[width= \textwidth]{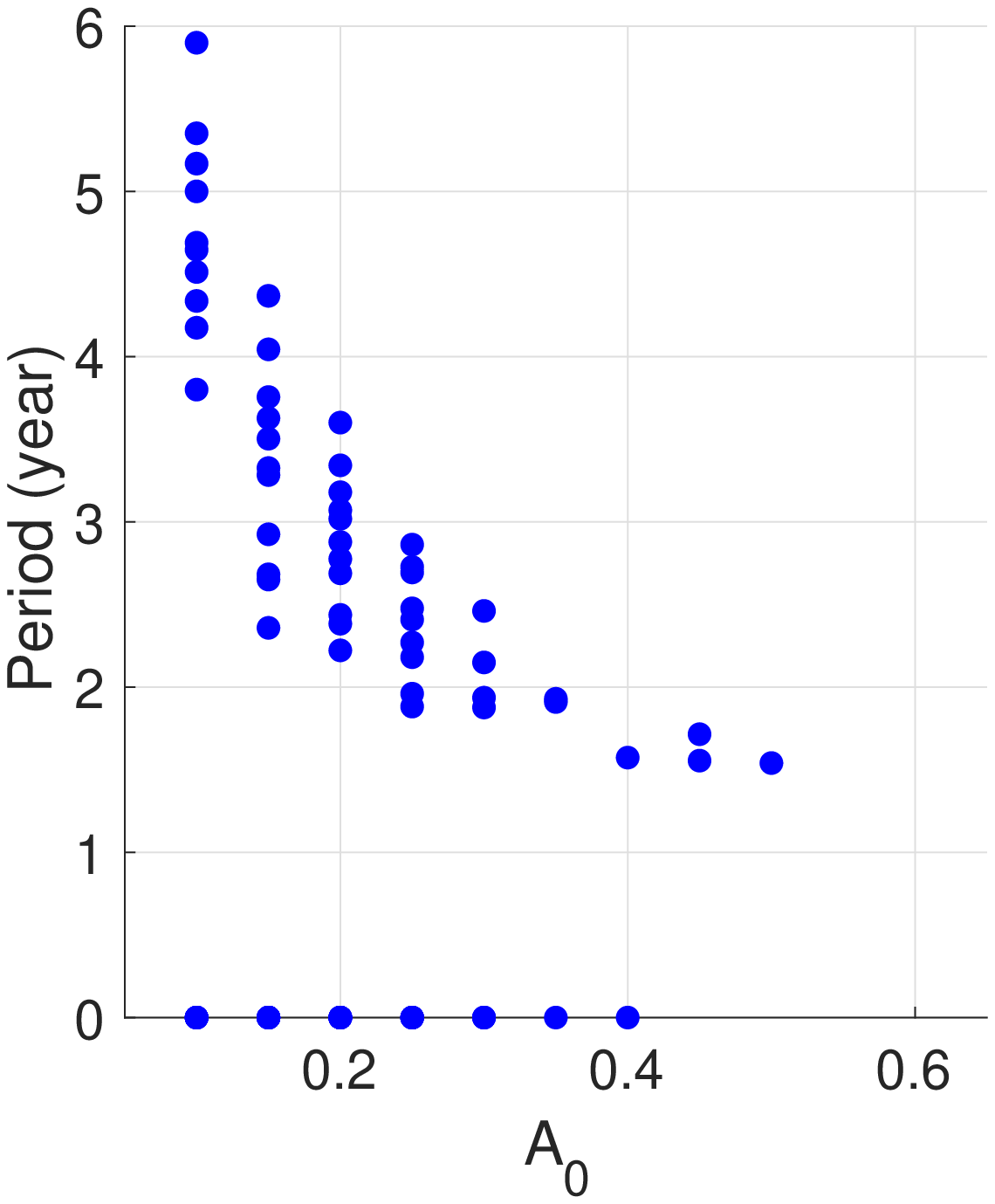}
		\caption{Changing $A_0$.}
	\end{subfigure}
	\begin{subfigure}{.32\textwidth}
		\centering
		\includegraphics[width= \textwidth]{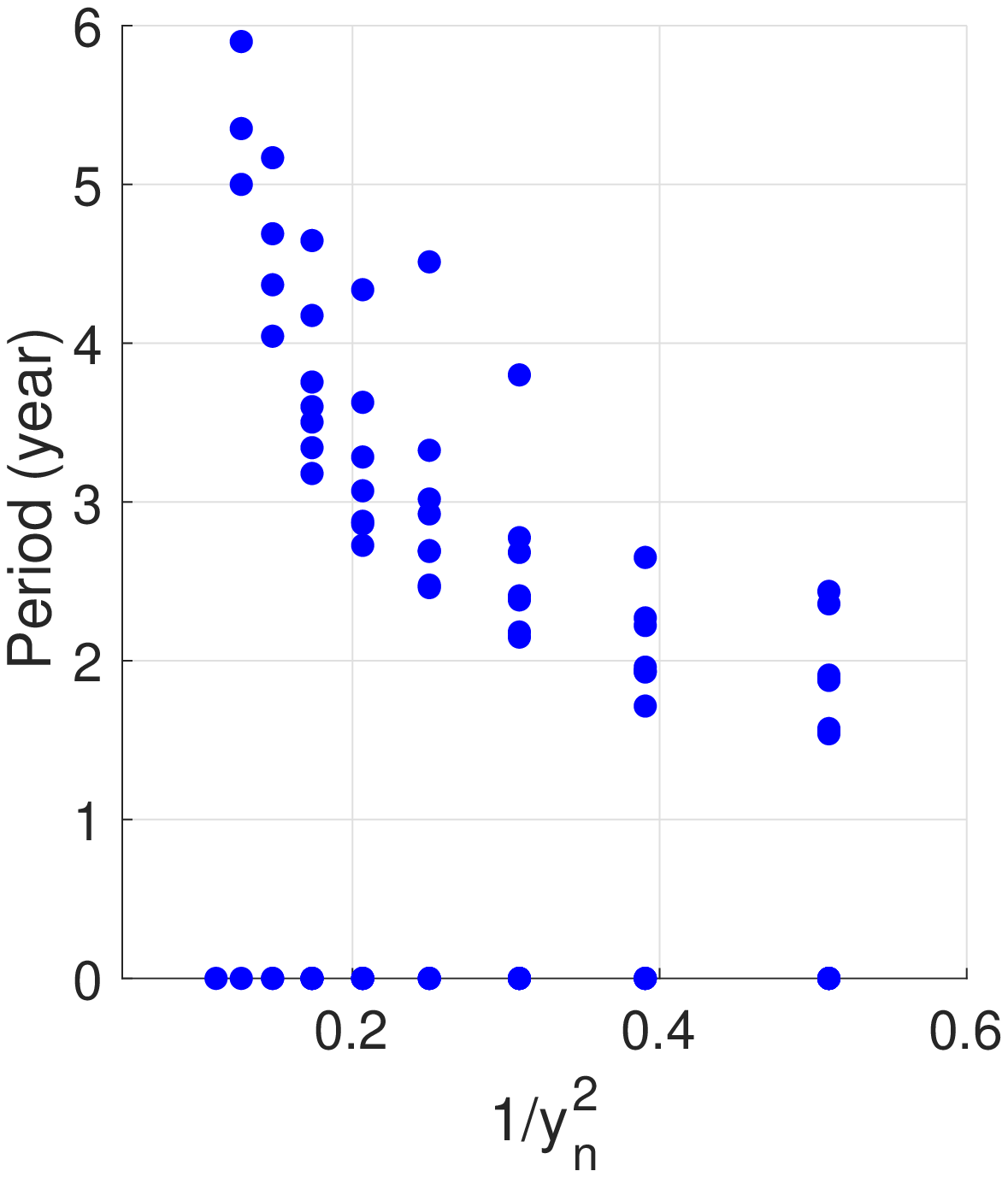}
		\caption{Changing $y_n$.}
	\end{subfigure}
	\caption{\label{fig:perdep}The dependence of the period of the oscillation on the parameters $\theta$, $A_0$ and $y_n$ for the VoC model. The range for each value of one of the parameters is due to the variation of the other two parameters.}
\end{figure}

We now turn to the scaling and nondimensionalization process to analyze the dependence on variable physical parameters in some detail. Only the VoC model is considered when comparing with the S\&S model, since the MZ model does not show periodic behaviour for most realistic parameter values.  The main uncertainties are in the latitude at which the Rossby waves travel ($y_n$), the wind forcing factor at $y_n$ ($\theta$), and the overall strength of the wind forcing ($A_0$). In Section \ref{ssec:bif} we used the values $y_n=2$, $A_0=0.2$, and $\theta=3$ for the reference scaling. In Figures \ref{fig:perch} and \ref{fig:perdep} these three parameters are varied according to Table \ref{tab:pardep} to study the dependence of the period. Note that these values differ from the those in Jin \cite{Jin1997a} due to a different scaling. All other parameter values are kept constant and listed in Appendix B.

\begin{table}[h]
	\centering
	\caption{\label{tab:pardep}The ranges in which the nondimensional parameters are varied. Values of $y_n=2$, $A_0=0.2$, and $\theta=3$ were used in Subsection \ref{ssec:bif}.}
	\begin{tabular}{llccc}
		\hline
		Parameter & & Dimensional & Dimensionless & Step\\
		\hline \hline
		Wind forcing factor at $y_n$ & $\theta$ & -  & 1.0 - 4.0 & 0.2 \\
		Wind forcing strength & $A_0$ & 0.5 - 3.0 $\cdot10^{-2}$ Pa & 0.1 - 0.6 & 0.05 \\
		Latitude Rossby waves & $y_n$ & $5.0^\circ$ - $12.1^\circ$ & 1.4 - 3.4 & 0.2 \\
		\hline
	\end{tabular}
\end{table}

Varying the parameters in the ranges given in Table \ref{tab:pardep}, the period of oscillation is computed for both the VoC model and the S\&S model. Note that the considered points can lie in the multistable region between the connecting orbit and Hopf curve. Figure \ref{fig:perch} compares the period of oscillation of both models. The period of the VoC model is larger whenever an oscillation is present. There are locations in parameter space where no oscillation occurs in the VoC model, where it does occur in the S\&S model. This can also be seen by comparing Figures \ref{fig:2param_SS} and \ref{fig:bifgamVoC}, where both the Hopf curve and connecting orbit curve have shifted upward for the VoC model. This upward shift reduces the region in parameter space in which periodic orbits occur.

The dependence of the period on the different parameters is shown in Figure \ref{fig:perdep}. The different points for each parameter value are due to the variation of the other two parameters. For increasing $\theta$, that is, when the effect of the wind forcing at higher latitudes increases, the period slowly increases. There turns out to be a minimum value for $\theta$ around 1.7 below which no oscillations occur. In that case the signal at higher latitudes is too weak to have a significant effect at the eastern boundary. Considering the strength of the wind forcing $A_0$, a stronger wind results in a shorter period. This could be due to the larger absolute difference between the effect at the equator and at higher latitudes, leading to a weaker effect of the latter. This decrease with increasing wind strength appears to be approximately exponential. For a realistic $A_0$ in the centre of the range, the period of the oscillation is approximately 2.5 to 3.5 years. This is still smaller than that of ENSO.

Looking at the latitude $y_n$ at which the Rossby waves travel, instead of the latitude itself, $1/y_n^2$ is plotted, since this gives the velocity of the Rossby wave traveling at that latitude. For higher velocities, so lower latitudes, the oscillations have a smaller period. The faster the wave travels, the shorter the delay is, resulting in a shorter period. Similarly, slow waves result in longer periods.

\section{Summary, Discussion and Conclusion}
\label{sec:disc}

Delay models are useful as conceptual climate models due to their infinite dimensional nature and reduced number of parameters. They are suited for mathematical analysis and therefore can add to the physical understanding of the processes involved. In this paper the Mori-Zwanzig formalism has been investigated as a method to derive delay equations. The two-strip model of the El Ni\~no Southern Oscillation (ENSO) was used as a test case for the application of the technique. The reason being that models incorporating delay have already been proposed for this phenomenon. 

The Mori-Zwanzig formalism gives an exact rewriting of a system of ordinary differential equations \cite{Chorin2002}. The rewritten equation contains a Markovian, noise, and memory term. Here the focus was on the memory term, since this is an integral over the history of the system, just as a delay term represents this history. The memory integral can be rewritten in the form of a distributed delay term by making several approximations. Under some further approximations this can be reduced to a discrete delay. The memory kernel determines the times at which the integral is large and this way gives the delay time. Because of the integral, a peak in the memory kernel does not necessarily coincide with the dominant delay time.

For the model of ENSO, in addition to the use of the Mori-Zwanzig formalism, also the method of variation of constants has been employed. Starting from the linear two-strip model \cite{Jin1997a}, the linear part of the delay model by Suarez and Schopf has been derived \cite{Suarez1988}. In this case the Mori-Zwanzig formalism is equivalent to using variation of constants (discussed in Section \ref{sec:mz}). Furthermore, a nonlinear version of the two-strip model was derived by assuming that the sea surface temperature is proportional to the subsurface temperature. From this nonlinear two-strip model (nonlinear) delay models have been derived. Here the two methods do not yield the same result. This is due to the approximations needed to obtain a closed-form equation from the Mori-Zwanzig formalism; whereas the method of variation of constants is exact.

Both derived nonlinear delay models contain an extra cubic delay term compared to the model by Suarez and Schopf \cite{Suarez1988}. In both cases this additional term results in an increased period of the model oscillation. Another consequence of the additional term is the decrease of the area in parameter space where stable oscillations occur. For the model derived using the Mori-Zwanzig formalism, this decrease is so large that no stable periodic behaviour occurs for most realistic values of the parameters. The exact derived model does show oscillations for these parameter values. The period of this model derived using variation of constants is closer to the real period of ENSO than the model proposed by Suarez and Schopf. However, its period is still smaller than what is seen in observations. One approach to improve the match between model period and data, could be to no longer assume a delta-function for the spatial pattern of the wind forcing, but rather take a more realistic pattern. As a consequence, the resulting delay model will no longer contain a discrete delay, but a distributed delay. Another option for possible improvement is to incorporate additional nonlinearities, for example in the thermocline equations. 
For the ENSO model it is not necessary to use the Mori-Zwanzig formalism to arrive at a delay equation. It does not give additional understanding compared to the method of variation of constants. This holds true for all models that are linear in the unresolved variables. This method yields a delay equation  for the resolved variable (temperature) in both the linear and nonlinear case, because the thermocline depth equations are linear in both models. The Mori-Zwanzig formalism gives the same result for the linear system, but not for nonlinear systems. There the formalism does not give accurate results, since approximations are needed.

When the equations of the model considered are also nonlinear in the unresolved variables, the Mori-Zwanzig formalism is the only method that will give a result. In such a nonlinear case the orthogonal dynamics system has to be approximated. This approximation needs to be an improvement on the pseudo-orthogonal dynamics approximation \cite{Gouasmi2017}, since this approximation was shown to be inaccurate for the ENSO model. It does not yield accurate results because the time scale of the unresolved variables is of the same order as that of the resolved variables. Only if its accuracy can be shown for a specific model, it can be justified to apply this approximation method. Additionally, a comment should be made on the computational cost of the pseudo-orthogonal dynamics approximation.  If the corresponding ODE system derived for the pseudo-orthogonal dynamics is of high dimension, then solving such a system may not prove to be an improvement over solving the original full ODE model. However, with regard to the justification of the use of DDE models for parameter studies, one is only interested in whether or not the resulting memory kernels have pronounced peaks. Hence, the pseudo-orthogonal dynamics approximation may only need to be solved once.  We expect this approach to be fruitful in climate models whenever the unresolved subsystem describes approximately linear wave transport. The derivation of improved approximations is a first step that needs to be taken to apply the Mori-Zwanzig formalism accurately to nonlinear models. This is a necessary step to be able to reliably derive nonlinear delay models for climate models.

Since many climate models are wave equations in one form or another, it is expected that they can also be represented by a delay system. This would imply that there is an abundance of phenomena in the climate system that can be described by a delay equation, possibly enabling the understanding of complex behaviour using relatively simple conceptual models. \\[3mm]




{\noindent
	\emph{Data access}: {The codes supporting this article have been uploaded as part of the supplementary material. They can also be found on the online repository figshare: https://doi.org/10.6084/m9.figshare.8085683.v1} \\[3mm]
%
\emph{Authors contributions}: {SF, CQ and HD designed the study. The work was carried out 
mainly by SF, with the exception of the bifurcation analysis in Section \ref{ssec:bif} which was done by CQ.  All authors contributed to the work, discussed the results and read and approved the 
manuscript.} \\[3mm]
\emph{Competing Interests}: {The author(s) declare that they have no competing interests
	.}\\[3mm]
%
\emph{Funding}: {This work was supported by funding from the European Union Horizon 2020 
research and innovation programme for the ITN CRITICS under Grant Agreement 
Number 643073 (CQ, JS and HD), the Mathematics of Planet Earth program (project No. 657.014.006) of the Dutch Science  Foundation (NWO, JF) and EPSRC Grants No. EP/N023544/1 and No. EP/N014391/1 (JS).}\\[3mm]
\emph{Acknowledgments}: {We thank the two anonymous reviewers for their constructive comments. SF thanks the University of  Exeter for hosting her for 5 months in 2018.}
}


\bibliographystyle{unsrt}

\begin{thebibliography}{10}

\bibitem{Battisti1988a}
D~S Battisti and A~C Hirst.
\newblock {Interannual Variability in a Tropical Atmosphere-Ocean Model:
  Influence of the Basic State, Ocean Geometry and Nonlinearity}.
\newblock {\em Journal of the Atmospheric Sciences}, 46(12):1687--1712, 1988.

\bibitem{Beck2009}
C~L Beck, S~Lall, T~Liang, and M~West.
\newblock {Model reduction, optimal prediction, and the Mori-Zwanzig
  representation of Markov chains}.
\newblock In {\em Joint 48th IEEE Conference on Decision and Control and 28th
  Chinese Control Conference}, 2009.

\bibitem{Chorin2002}
A~J Chorin, O~H Hald, and R~Kupferman.
\newblock {Optimal prediction with memory}.
\newblock {\em Physica D}, (166), 2002.

\bibitem{Darve2009}
E~Darve, J~Solomon, and A~Kia.
\newblock {Computing generalized Langevin equations and generalized
  Fokker-Planck equations}.
\newblock {\em Proceedings of the National Academy of Sciences},
  106(27):10884--10889, 2009.

\bibitem{DijkstraBook}
Henk~A Dijkstra.
\newblock {\em {Nonlinear Physical Oceanography}}, volume~28.
\newblock Springer, 2nd revise edition, 2005.

\bibitem{Dijkstra1995}
Henk~A Dijkstra and J~D Neelin.
\newblock {On the attractors of an intermediate coupled ocean-atmosphere
  model}.
\newblock {\em Dynamics of Atmosphere and Ocean}, (19-48), 1995.

\bibitem{Dominy2017}
Jason~M. Dominy and Daniele Venturi.
\newblock {Duality and conditional expectations in the Nakajima-Mori-Zwanzig
  formulation}.
\newblock {\em Journal of Mathematical Physics}, 58(8), 2017.

\bibitem{DDEbiftool1}
Koen Engelborghs, Tatyana Luzyanina, and Dirk Roose.
\newblock {Numerical bifurcation analysis of delay differential equations using
  DDE-BIFTOOL}.
\newblock {\em ACM Transactions on Mathematical Software (TOMS)}, 28(1):1--21,
  2002.

\bibitem{DDEbiftool2}
Koen Engelborghs, Tatyana Luzyanina, and Giovanni Samaey.
\newblock {DDE-BIFTOOL v. 2.00: a Matlab package for bifurcation analysis of
  delay differential equations}.
\newblock 2001.

\bibitem{PDEBook}
L~C Evans.
\newblock {\em {Partial Differential Equations}}, volume~19.
\newblock American Mathematical Society, 2nd edition, 2010.

\bibitem{Givon2004}
D~Givon, R~Kupferman, and A~Stuart.
\newblock {Extracting macroscopic dynamics: model problems and algorithms}.
\newblock {\em Nonlinearity}, 17:55--127, 2004.

\bibitem{Gottwald2017}
Georg~A Gottwald, Daan~T Crommelin, and Christian L~E Franzke.
\newblock {\em {Stochastic Climate Theory}}, pages 209--240.
\newblock Cambridge University Press, 2017.

\bibitem{Gouasmi2017}
A~Gouasmi, E~J Parish, and K~Duraisamy.
\newblock {A priori estimation of memory effects in reduced-order models of
  nonlinear systems using the Mori-Zwanzig formalism}.
\newblock {\em Proceedings of the Royal Society A}, (20170385), 2017.

\bibitem{Hao1993}
Z~Hao, J~D Neelin, and F~Jin.
\newblock {Nonlinear Air-Sea Interaction in the Fast-Wave Limit}.
\newblock {\em Journal of Climate}, 6:1523--1544, 1993.

\bibitem{Jin1997a}
F~Jin.
\newblock {An Equatorial Ocean Recharge Paradigm for ENSO. Part I: Conceptual
  Model}.
\newblock {\em Journal of Atmospheric Sciences}, 54:811--829, 1997.

\bibitem{Jin1997b}
F~Jin.
\newblock {An Equatorial Ocean Recharge Paradigm for ENSO. Part II: A
  Stripped-Down Coupled Model}.
\newblock {\em Journal of Atmospheric Sciences}, 54:830--847, 1997.

\bibitem{Kato1995}
T.~Kato.
\newblock {\em {Perturbation Theory for Linear Operators}}.
\newblock Springer-Verlag, 1995.

\bibitem{Keane2017}
A~Keane, B~Krauskopf, and C~M Postlethwaite.
\newblock {Climate models with delay differential equations}.
\newblock {\em Chaos}, 27(114309), 2017.

\bibitem{Krauskopf2014}
B~Krauskopf and J~Sieber.
\newblock {Bifurcation analysis of delay-induced resonances of the
  El-Ni{\~{n}}o Southern Oscillation}.
\newblock {\em Proceedings of the Royal Society A}, 470(2169), 2014.

\bibitem{Mori1965}
H~Mori.
\newblock {Transport, Collective Motion and Brownian Motion}.
\newblock {\em Progress of Theoretical Physics}, 33(3):423--455, 1965.

\bibitem{morriss2013}
Gary~P Morriss and Denis~J Evans.
\newblock {\em {Statistical Mechanics of Nonequilbrium Liquids}}.
\newblock ANU Press, 2013.

\bibitem{Runge2014}
J~Runge, V~Petoukhov, and J~Kurths.
\newblock {Quantifying the Strength and Delay of Climate Interactions: The
  Ambiguities of Cross Correlation and Novel Measure Based on Graphical
  Models}.
\newblock {\em Journal of Climate}, 27:720--739, 2014.

\bibitem{Shen2011}
J~Shen, T~Tang, and L~Wang.
\newblock {\em {Spectral Methods}}.
\newblock Springer.

\bibitem{DDEbiftool3}
Jan Sieber, Koen Engelborghs, Tatyana Luzyanina, Giovanni Samaey, and Dirk
  Roose.
\newblock {DDE-BIFTOOL Manual-Bifurcation analysis of delay differential
  equations}.
\newblock {\em arXiv preprint arXiv:1406.7144}, 2014.

\bibitem{Suarez1988}
M~J Suarez and P~S Schopf.
\newblock {A Delayed Action Oscillator for ENSO}.
\newblock {\em Journal of the Atmospheric Sciences}, 45(21):3283--3287, 1988.

\bibitem{Szalai2014}
R~Szalai.
\newblock {Modelling elastic structures with strong nonlinearities with
  application to stick-slip friction}.
\newblock {\em Proceedings of the Royal Society A}, (470), 2014.

\bibitem{Tziperman1998}
Eli Tziperman, Mark~A. Cane, Stephen~E. Zebiak, Yan Xue, and B.~Blumenthal.
\newblock {Locking of El Nino's peak time to the end of the calendar year in
  the delayed oscillator picture of ENSO}.
\newblock {\em Journal of Climate}, 11(9):2191--2199, 1998.

\bibitem{Wouters2016}
J~Wouters, S~I Dolaptchiev, V~Lucarini, and U~Achatz.
\newblock {Parametrization of stochastic multiscale triads}.
\newblock {\em Nonlinear Processes in Geophysics}, 23:435--445, 2016.

\bibitem{Zebiak1987}
S~E Zebiak and M~A Cane.
\newblock {A Model of El Ni{\~{n}}o-Southern Oscillation}.
\newblock {\em Monthly Weather Review}, 115:2262--2278, 1987.

\bibitem{Zhu2018a}
Yuanran Zhu, Jason~M. Dominy, and Daniele Venturi.
\newblock {On the estimation of the Mori-Zwanzig memory integral}.
\newblock {\em Journal of Mathematical Physics}, 59(10), 2018.

\bibitem{Zhu2018b}
Yuanran Zhu and Daniele Venturi.
\newblock {Faber approximation of the Mori–Zwanzig equation}.
\newblock {\em Journal of Computational Physics}, 372:694--718, 2018.

\bibitem{Zwanzig1973}
R~Zwanzig.
\newblock {Nonlinear Generalized Langevin Equations}.
\newblock {\em Journal of Statistical Physics}, 9(3):215--220, 1973.

\end{thebibliography}

\clearpage

\appendix

\section{Pseudo-Orthogonal Dynamics Approximation}
\label{app:pod}

\subsection{Formulation}
\label{ssec:pod_form}

To simplify the issue of solving the orthogonal dynamics system, Gouasmi \emph{et al.} \cite{Gouasmi2017} derived the pseudo-orthogonal dynamics (POD) equation. Under certain assumptions this is an exact rewriting of the orthogonal dynamics equation. In this rewritten form the orthogonal dynamics system can be more easily solved. The main assumption in the approach is the commutativity of $\e^{tQ\mathcal{L}}$ and $R$:
\begin{equation}
[\e^{tQ\mathcal{L}}R](x) \approx [R\e^{tQ\mathcal{L}}] (x) = R(\phi^Q(x,t)).
\end{equation}
For linear systems commutativity holds and this relation is exact, as shown in \cite{Gouasmi2017}. For nonlinear systems the above relation may be used as an approximation. The accuracy of this approximation is not a priori clear and requires verification. With the assumption of commutativity, the orthogonal dynamics equation can be reformulated into the POD equation:
\begin{equation}
\label{eq:pseudo}
\frac{\partial}{\partial t}\phi^Q(x,t) = R(\phi^Q(x,t)) - R(\hat{\phi}^Q(x,t)).
\end{equation}
Note that this equation can be implemented more straightforward in numerical codes compared to the original orthogonal dynamics equation. When this equation is solved, the noise term corresponding to the resolved component $\phi_i$ is given by
\begin{equation}
\label{eq:noise_pseudo}
F_i(x_0,t) = R_i(\phi^Q(x,t)) - R_i(\hat{\phi}^Q(x,t)).
\end{equation}
This is the part of the POD system corresponding to the respective resolved variable. The noise term thus can be retrieved directly when solving the POD equation. Defining $R_Q(x) = R(x) - [PR](x)$, the error made in this approximation is:
\begin{equation}
\delta \leq \big|  [\partial_x R_{Qi}](x) - [\partial_x R_{Qi}](\hat{x}) \big|
\end{equation}

The POD approximation simplifies solving the orthogonal dynamics system. However having a solution to this system is not necessarily sufficient to also formulate an expression for the memory term. When the noise term is a complicated function of the solution to the orthogonal dynamics equation, it can be quite difficult to get an expression for the memory kernel. Therefore, it is useful to look into approximations of the kernel. The following approximation has been derived by Gouasmi \emph{et al.} \cite{Gouasmi2017}. The first step is to consider the $n$ components of $\mathcal{L}F_i(x,t)$ as the partial derivative of $F_i(x,t)$ in the direction of $\bar{R}(x) = R(x)/||R(x)||$, instead of $n$ separate derivatives in directions $R_j(x)$ for each coordinate $x_j$ with $j=1,\dots,n$. Here $\bar{R}:\mathbb{R}^n\rightarrow\mathbb{R}^n$ and $||\cdot||$ is the $l^2$-norm, i.e.~the standard Euclidean norm on $\mathbb{R}^n$. This yields
\begin{equation}\label{eq:memoryapprox}
\mathcal{L}F_i(x,t) = \sum_{j=1}^N R_j(x) \partial_{x_j}F_i(x,t) = ||R(x)|| \lim_{\epsilon\rightarrow0}\frac{F_i(x+\epsilon\bar{R}(x),t)-F_i(x,t)}{\epsilon}.
\end{equation}
Here the right-hand side is written as the limit corresponding to the derivative of $F_i$ at $x$ in the direction of $\bar{R}(x)$. The memory kernel $P\mathcal{L}F_i$ is the projection of this equation onto the resolved variables. Since $F_i(x,t)$ is the solution to the orthogonal dynamics equation, it only depends on the unresolved variables. Therefore, the second term in the numerator of (\ref{eq:memoryapprox}) disappears after projection onto the resolved variables, i.e.~$PF_i(x,t) = F_i(\hat{x},t)=0$. Applying $P$ to Equation (\ref{eq:memoryapprox}) and using finite differences, results in the following approximation for the memory kernel:
\begin{equation}
\label{eq:kernelapprox}
K_i(\hat{x},t) \approx ||R(\hat{x})|| \frac{F_i(\hat{x}+\epsilon\bar{R}(\hat{x}),t)}{\epsilon}.
\end{equation}
The exact result can be recovered if the limit $\epsilon\rightarrow0$ exists. The memory integral can be approximated by the rectangle rule or another approximation method for integrals. 

The derived approximation is most useful when working numerically. To obtain analytical expressions using this method can be cumbersome. However, if the result of the approximation remains tractable and the limit can be computed, it can result in an exact expression for the memory kernel.

Gouasmi \emph{et al.} \cite{Gouasmi2017} applied this approximation to numerically compute the memory for the POD system to the Burger's equation and the Kuramoto-Sivashinsky equation. They found that the approximation in the Burger's equation was accurate for high wave numbers, but not for low wave numbers. This means it gives a good approximation when the memory is determined mainly by fast dynamics and thus is relatively short. However, for systems where slow dynamics affects the memory term, resulting in a longer memory, it appears to be less accurate.

\subsection{Application to Nonlinear ENSO Model}
\label{ssec:pod_ENSO}

Since it is not feasible to analytically solve the orthogonal dynamics equation for the nonlinear ENSO model as discussed in Section \ref{ssec:nonlts}, we apply the POD approximation to find an approximate solution. The conditions for this approximation are not met for the nonlinear two-strip model, but it can be used as a first estimate.

The POD equations for the nonlinear two-strip model of equation \eqref{eq:twostr} are
\begin{equation}
\label{eq:mznlpseu}
\begin{split}
\partial_t h_c^Q(x,t) &= -(\epsilon_0+\partial_x)h_c^Q(x,t), \\
\partial_t h_n^Q(x,t) &= -(\epsilon_0-\frac{1}{y_n^2}\partial_x)h_n^Q(x,t), \\
\partial_t T_e^Q(x,t) &= c_{h}^*(x)\big(1 - \beta T_e^{Q}(x,t)^2\big) \big(h_c^Q(x,t) + \frac{1}{1+y_n^2} h_n^Q(x,t)\big).
\end{split}
\end{equation}
The first two equations have exponential functions as solutions. To find a solution for $T_e^Q$ we substitute the solutions for $h_c^Q$ and $h_n^Q$ into the equation for $T_e^Q$. The solution for $T_e^Q$ with the condition that $(T_e^{Q})^{2}<\frac{1}{\beta}$ is
\begin{equation}
\label{eq:TeQtanh}
\begin{split}
T_e^Q(x,t) &= \frac{1}{\beta}\tanh^2\Big(\arctanh\big(\sqrt{\beta}T_e(x,0)\big) \\
&\qquad + c_h^*(x)\sqrt{\beta}\big( (1-\e^{-(\epsilon_0 + \partial_x)t}) (\epsilon_0+\partial_x)^{-1} h_c(x,0)  \\
& \qquad + \frac{1}{1+y_n^2}(1-\e^{-(\epsilon_0- \frac{1}{y_n^2}\partial_x)t}) (\epsilon_0-\frac{1}{y_n^2}\partial_x)^{-1} h_n(x,0) \big) \Big).
\end{split}
\end{equation}
If $(T_e^{Q})^{2}>\frac{1}{\beta}$ the $\tanh$ has to be replaced by a $\coth$, and when $(T_e^{Q})^{2}=\frac{1}{\beta}$ the result is a constant $T_e^Q$, since then $\partial_t T_e^Q(x,t)=0$. The initial conditions determine which of the solutions should be used. Most likely is that $(T_e^{Q})^{2}<\frac{1}{\beta}$, as $\beta$ is small and $T_e$ is of order one. Therefore, in the following Equation (\ref{eq:TeQtanh}) is used.

The noise term is given by the right-hand side of the equation for $T_e^Q$ in Equation (\ref{eq:mznlpseu}), for which now a closed expression is known. The next step is to compute the memory kernel. By the presence of the hyperbolic tangent and several nonlinearities, Equation (\ref{eq:kernelapprox}) is used to approximate the memory kernel. Taking the limit $\epsilon\rightarrow 0$, yields the memory kernel. The equation for $T_e$ in the POD approximation becomes
\begin{equation}
\begin{split}
\frac{\d T_e}{\d t}(x,t) &= -c_T(x) T_e(x,t) + c_{h}^*(x)\big(\e^{-(\epsilon_0 + \partial_x)t}h_c(x,0)+\frac{1}{1+y_n^2}\e^{-(\epsilon_0 - \frac{1}{y_n^2}\partial_x)t}h_n(x,0)\big) \\
& \qquad \cdot \big(1 - \beta T_e^{Q}(x,t)^2\big) \\
& \quad + \int_0^t c_{h}^*(x) \big( 1+ \beta T_e^2(x,s) \big) \big( \mu \big( 1- \frac{\theta}{1+y_n^2} \big) \e^{-(\epsilon_0 + \partial_x)(t-s)}  g(x) T_e(x_E,s) \\
& \qquad - \mu \frac{\theta}{y_n^2}\frac{A_{rW}}{1+y_n^2}\e^{-(\epsilon_0 - \frac{1}{y_n^2}\partial_x)(t-s)}   g(x) T_e(x_E,s)  \big)  \d s.
\end{split}
\end{equation}

Similar to the procedure followed for the linear model in Sections \ref{ssec:char} and \ref{ssec:dellin}, this equation can be simplified. The desired result is an equation for the temperature in the east of the basin. Applying the method of characteristics one gets rid of the exponential $\partial_x$-terms. Assuming no reflection takes place at the eastern boundary, the noise term vanishes. The memory term further simplifies to two delay terms when a localized wind forcing is assumed. Considering the short delay as being instantaneous, the resulting nonlinear delay equation is
\begin{equation}
\begin{split}
\frac{\d T_e^E}{\d t} &= (c_S^*-c_T(x_E)) T_e^E(t) - c_L^* T_e^E(t-d) - \beta c_S^* T_e^{E}(t)^3 + \beta c_L^* T_e^{E}(t-d)^3,
\end{split}
\end{equation}
where
\begin{equation}
\begin{split}
c_S^* &= \mu A_0 \big(1-\frac{\theta}{1+y_n^2}\big) c_h^*(x_E) \e^{-\epsilon_0 (1-x_w)}, \\
c_L^* &= \mu A_0 \frac{\theta}{y_n^2}\frac{A_{rW}}{1+y_n^2} c_h^*(x_E) \e^{-\epsilon_0 (1+ y_n^2 x_w)}, \\
d &= 1+y_n^2x_w.
\end{split}
\end{equation}

\section{Scaling ENSO Model}
\label{app:scaling}

In this appendix we give an overview of the parameters and scaling used in the final delay model of Equation \eqref{eq:delscvoc}, which is repeated here:
\begin{equation}
\frac{\d \tilde{T}}{\d \tilde{t}} = \tilde{T}(\tilde{t})-\tilde{T}(\tilde{t})^3 - \alpha \tilde{T}(\tilde{t}-\delta) \big( 1 - \gamma \tilde{T}(\tilde{t})^2 \big),
\end{equation}
where 
\begin{equation}
\alpha = \frac{c_L^*}{c_S^*-c_T(x_E)}, \qquad \gamma = \frac{c_S^*-c_T(x_E)}{c_S^*}, \qquad \delta = (c_S^*-c_T(x_e))d.
\end{equation}
To arrive at this equation time and temperature are  scaled by:
\begin{equation}
\tilde{t} = (c_S^*-c_T(x_E)) t, \qquad \tilde{T} = \sqrt{\frac{\beta c_S^*}{c_S^*-c_T(x_E)}}T.
\end{equation}
Note that in Section \ref{ssec:delan} we drop the tildes after scaling.

The parameters defining $\alpha$, $\gamma$ and $\delta$ are:
\begin{equation}
\begin{split}
c_S^* &= \mu A_0 \big(1-\frac{\theta}{1+y_n^2}\big) c_h^*(x_E) \e^{-\epsilon_0 (1-x_w)}, \\
c_L^* &= \mu A_0 \frac{\theta}{y_n^2}\frac{A_{rW}}{1+y_n^2} c_h^*(x_E) \e^{-\epsilon_0 (1+ y_n^2 x_w)}, \\
d &= 1+y_n^2x_w,
\end{split}
\end{equation}
where
\begin{equation}
\begin{split}
c_T(x) &= \epsilon_w + 0.5\Big(1-\alpha_0 + (1+\alpha_0)\tanh\big(\frac{\delta_F^1}{\epsilon}F(x)\big)\Big) \delta_F^1F(x), \\
c_{h}^*(x) &= 0.5 \big(\tanh(\frac{\delta_F^1}{\epsilon}F(x)) - 1 \big) \alpha_0 \delta_F^1 F(x) (T_{0}-T_{s0})\frac{H}{H^*},
\end{split}
\end{equation}
with background wind forcing
\begin{equation*}
F(x)=0.6\Big( 0.12 - \cos\big(\frac{x-x_0}{2x_0}\pi\big)^2 \Big),
\end{equation*}
for $x_0=0.57$ and parameters 
\begin{equation}
\epsilon_w = \frac{\epsilon_T L}{c_0}, \qquad \alpha_0=\frac{H_1}{\tilde{H}}, \qquad \delta_F^1 = \frac{\tau_0 L}{c_0} \frac{b_w}{H_1}.
\end{equation}
The values of the involved parameters above are given in Table \ref{tab:val}. Here the eastern boundary is not considered at $x=1$, but at $x=x_E$ to avoid boundary effects of the model. The other dimensionless variables in the definitions of $c_S^*$ and $c_L^*$ are $\mu$, $\epsilon_0$, $A_{rW}$, $\theta$, $y_n$, $A_0$ and $x_w$. Based on the book \textit{Nonlinear Physical Oceanography} by Dijkstra \cite{DijkstraBook}
, we set
\begin{equation}
\mu=1, \qquad \epsilon_0 = \frac{a_M L}{c_0}, \qquad A_{rW} = r_W(1+y_n^2) - 1, \qquad x_w = 0.6,
\end{equation}
with $r_W=3/5$. The other parameters are set as $\theta=3$, $y_n=2$, and $A_0=0.2$ in Section \ref{ssec:bif} and varied as discussed in Section \ref{ssec:pardep} to study the dependence of the period.

\setcounter{table}{1}
\begin{table}[h]
	\centering
	\caption{\label{tab:val}Parameter values used in determining $c_T$ and $c_h^*$.}
	\begin{tabular}{lll}
		\hline
		Damping scale Newtonian cooling & $\epsilon_T$ & $9.25\cdot 10^{-8}$ s$^{-1}$ \\
		Basin length & $L$ & $1.5\cdot 10^7$ m \\
		Velocity first baroclinic Kelvin mode & $c_0$ & 2 m/s \\
		Background wind forcing strength & $\tau_0$ & $2.667\cdot 10^{-7}$ m/s$^2$ \\
		Parametrization constant & $b_w$ & $1.026\cdot10^2$ s \\
		Depth surface layer & $H_1$ & 50 m \\
		Depth top layer & $H$ & 200 m \\
		Depth for temperature gradient & $\tilde{H}$ & 50 m \\
		Steepness transition subsurface temperature & $H^*$ & 30 m \\
		Temperature without dynamics & $T_0$ & 30 $^\circ$C \\
		Background subsurface temperature & $T_{s0}$ & 22 $^\circ$C \\
		Rayleigh friction coefficient & $a_M$ & $1.3\cdot 10^{-8}$ s$^{-1}$\\
		Scaling parameter & $\epsilon$ & $10^{-4}$ \\
		Reference point in east of basin & $x_E$ & 0.9 \\
		\hline
	\end{tabular}
\end{table}

We also note the nondimensionalisation of the shallow water equations used in \cite{DijkstraBook}:
\begin{equation}
\tilde{y}  =  \sqrt{\frac{c_0}{\beta}} y, \quad \tilde{x} = L x, \quad \tilde{t} = \frac{L}{c_0} t, \quad \tilde{h} = H h, \quad \tilde{u} = c_0 u, \quad \tilde{v} = c_0 v \quad \text{and} \quad \tilde{\tau} = \frac{\rho H c_0^2}{L} \tau
\end{equation}
where $\rho$ the density, $\beta$ the beta-plane parameter and the tildes denote the dimensional quantities. This scaling is also needed for the dimensionalisation of the parameter values in Table \ref{tab:pardep}.

\end{document}